\documentclass[review, font=11]{elsarticle}

\usepackage{lineno,hyperref}
\modulolinenumbers[5]

\journal{Preprint submitted to Applied Energy}

\graphicspath{{pics/}}
\usepackage[cmex10]{amsmath}
\usepackage{amssymb,amsthm,accsupp}
\usepackage[linesnumbered,ruled,vlined]{algorithm2e}
\usepackage{stmaryrd}
\usepackage{fixltx2e}
\usepackage{color}

\usepackage{multirow,booktabs}
\usepackage{algpseudocode}
\usepackage[linesnumbered,ruled,vlined]{algorithm2e}
\usepackage{algcompatible}
\algblockdefx{ForP}{EndP}[1]%
  {\textbf{for }#1 \textbf{do in parallel}}%
  {\textbf{end for}}
\algblockdefx{ForPA}{EndPA}[1]%
  {\textbf{for all }#1 \textbf{do in parallel}}%
  {\textbf{end for}}

\usepackage[caption=false,font=footnotesize,labelformat=simple]{subfig}

\usepackage{textcomp}
\newcommand\figref{Fig.~\ref}

\DeclareMathOperator*{\argmin}{argmin}
\DeclareMathOperator*{\argmax}{argmax}
\newcommand\norm[1]{\left\lVert#1\right\rVert}

\newcommand{\RNum}[1]{\uppercase\expandafter{\romannumeral #1\relax}}
\usepackage{tabularx}
\makeatletter
\def\hlinewd#1{%
\noalign{\ifnum0=`}\fi\hrule \@height #1 %
\futurelet\reserved@a\@xhline}
\makeatother

\usepackage{threeparttable}

\begin{document}

\begin{frontmatter}

\title{Hierarchical Distributed Framework for EV Charging Scheduling Using Exchange Problem}

\author[mymainaddress]{Behnam Khaki\corref{mycorrespondingauthor}} 
\cortext[mycorrespondingauthor]{Corresponding author}
\ead{behnamkhaki@ucla.edu}

\author[mymainaddress]{Chicheng Chu}
\ead{peterchu@ucla.edu}

\author[mymainaddress]{Rajit Gadh}
\ead{gadh@ucla.edu}

\address[mymainaddress]{Smart Grid Energy Research Center (SMERC), University of California, Los Angeles\\
44-120 Engr.IV, 420 Westwood Plaza, Los Angeles, CA, USA, 90095}

\begin{abstract}
In this paper, a distributed trilayer multi-agent framework is proposed for optimal electric vehicle charging scheduling (EVCS). 
The framework reduces the negative effects of electric vehicle charging demand on the electrical grids.
To solve the scheduling problem, a novel hierarchical distributed EV charging scheduling (HDEVCS) is developed as the \textit{exchange problem}, where the agents are clustered based on their coupling constraints.
According to the separability of the agents' objectives and the clusters' coupled constraints, HDEVCS is solved efficiently in a distributed manner by the alternating direction method of multipliers. Comparing to the exiting trilayer methods, HDEVCS reduces the convergence time and the iteration numbers since its structure allows the agents to update their primal optimization variable simultaneously.
The performance of HDEVCS is evaluated by numerical simulation of two small- and large- scale case studies consisting of $306$ and $9051$ agents, respectively.
The results verify the scalability and efficiency of the proposed method, as it reduces the convergence time and iteration numbers by $60\%$ compared to the state-of-the-art methods, flattens the load profile and decreases the charging cost considerably without violating the grid feeders' capacity.
The significant outcome of our method is the accommodation of a large EV population without investment in grid expansion.
\end{abstract}


\end{frontmatter}


\section{Introduction}
\paragraph{Motivation and Problem Statement}
In the transportation sector, the transition from non-electric to electric-powered end-use consumers - \textit{Electrification}- is mostly influenced by the electric vehicles (EVs) (\citet{mai2018electrification}).
EVs significantly reduce greenhouse gas emissions in the transportation sector.
However, they may increase the peak load demand and energy loss in the residential buildings and distribution power grids (\citet{{Fischer1}, {Clement-Nyns1}, {Fernandez1}}) which result in inefficient operation of the system. Through EV load demand management, however, the distribution network operator (DNO) can not only decrease the negative effects of EV charging load without investment in grid capacity expansion, but also benefit by peak load shaving or load variance minimization (LVM). EVs, owing to the bidirectional EV chargers, can also play the role of energy storage and feed electricity back to the grid to provide regulation services. This is called vehicle-to-grid (V2G) which its potential was investigated by \citet{Kempton1}. Nevertheless, EV load management is challenging due to the uncertainty in arrival time, departure time and energy demand (\citet{8585744, 8440360}), limited capacity of the energy resources and distribution grid equipment (\citet{8440531}), the scalability issue, and EV owners' privacy preserving. In this paper, we propose an EV charging scheduling (EVCS) framework to address the scalability and privacy preserving issues.  

\paragraph{Related Work}
There is a rich body of literature proposing a variety of approaches for EVCS which fall into two categories: centralized and distributed. In the former, a central entity such as DNO receives all the required data over the communication system from the dispersed EVs and coordinates their charging demand. Then, the optimal charging profiles are sent back to EVs: By \citet{Clement-Nyns1}, it is shown that the uncoordinated EV charging increases power loss and voltage deviation significantly, therefore the authors propose a centralized method where the EV owners have no control over the charging profile, and it is decided by DNO; a model predictive based algorithm is proposed by \citet{Tang1} for total charging cost reduction, where the authors use the truncated sample average approximation to reduce the complexity of their centralized method at the cost of performance degradation; \citet{Bin1} introduce a centralized event-triggered receding horizon method to reduce EV charging cost in a campus parking; an optimal strategy for V2G aggregator is designed by \citet{Peng1} to maximize the economic benefit by participation in frequency regulation while satisfying EV owners' demand; a centralized algorithm is designed by \citet{Bilh1} to flatten the netalod fluctuations due to renewable energy resources using EV charging control; \citet{Zheng1} propose a real-time EV charging scheduling where the computational complexity is reduced by introducing a capacity margin and the charging priority indices; a centralized mechanism is proposed by \citet{Diaz1} in which a third-party entity coordinates a day-ahead bidding system to optimize the global bid; a transactive EV charging management is presented in \citet{Liu1} to maximize the real-time profit based on the net electricity exchange with the grid; and a two-layer centralized EVCS is proposed by \citet{Mehta1} where each aggregator optimizes active power of the EVs in the first layer, and the second layer provides reactive power management for loss reduction in the grid. The main issues with the centralized approaches are: (\RNum{1}) EV owners' privacy as they have to communicate their sensitive charging information (e.g. arrival and departure times as well as their battery energy level and capacity) with the central coordinator, (\RNum{2}) curse of dimensionality when EV penetration increases in the grid, and (\RNum{3}) vulnerability to a single point failure, i.e. if the central coordinator faces any problem, the whole system breaks down. In the distributed approaches, which address the above issues, the central entity coordinates the EV charging demand through communication with the EV charger agents or EV aggregators (EVAs). Instead of solving the scheduling problem centrally, it is solved through a distributed and iterative procedure in collaboration with the EVA and EV agents.

Among the distributed methods, two different EV charging infrastructures (ECIs) are considered: (\RNum{1}) bilayer structure which consists of either EVA and EVs or DNO and EVAs, and (\RNum{2}) trilayer structure which includes DNO, EVAs, and EVs. \citet{Low1} propose a distributed charging scheduling which is solved using the projected gradient descent and provides valley filling. In \citet{Rivera1}, the EVCS problem is formulated as the \textit{exchange problem} which is efficiently solved by the alternating direction method of multipliers (ADMM) for the valley filling and charging cost reduction (CR). \citet{Xiong1} use a water filling algorithm incentivizing the EV owners to shift their charging demand to the off-peak hours. A mean-field game theory-based method is proposed by \citet{Kebriaei1} to provide valley filling and reduce battery degradation cost. The work is further expanded by \citet{Kebriaei2} to consider plug-in hybrid EVs, in which the propulsion provided by gasoline gives more flexibility to EVs for participation in V2G. \citet{Latifi1} propose a distributed method where each EV only communicates with its neighbors to reduce the communication overhead, and the optimal charging problem is addressed by the full Nash Folk theorem. \citet{Ledwich1} propose a bi-level programming based hierarchical decomposition EVCS where the objective is to reduce the generation unit cost for DNO. \citet{Gupta1} propose three algorithms based on the projected gradient descent (\citet{Low1}) and ADMM, where EVAs communicate with their neighbors, and each EVA centrally calculates charging profiles of EVs which it supplies. The proposed distributed bilayer EVCS by \citet{Crespi1} reduces charging cost and increases EV owners' convenience, where EVAs update their charging load sequentially, not in parallel. \citet{Staudt1} propose a decentralized method where the EVAs receive signals from DNO if their aggregated EV load results in the transmission line congestion. EVAs help congestion relief by rescheduling all the EVs' load demands and available V2G capacity. Considering the on-site uncertain wind generation, a distributed method is proposed by \citet{Yang1} to increase the local wind energy utilization and satisfy the EV load demand through penalizing energy purchase from the grid.

Considering the trilayer EV charging structure, \citet{Ting1} use ADMM to schedule the EV charging demand with the purpose of user convenience maximization, which is characterized by the EV battery's final state of charge. A distributed approach based on the sub-gradient method is proposed by \citet{Song1} to satisfy customers' charging demands and the coupled constraints relating to the EVAs' feeder capacity. The EVA constraints, however, are relaxed in \cite{Song1} and included in the Lagrangian of the cost function which results in suboptimality of the coordinated charging profiles. \citet{Li1} design a distributed framework for V2G scheduling where EVAs' revenue is neglected, and the gradient projection method is used to solve optimal V2G scheduling problem which is not computationally efficient, and its key convergence parameter depends on the number of EVs. The hierarchical framework designed by \citet{Shao1} includes two iterative procedures: the first procedure is between DNO and EVAs, and the second one is between each EVA and its EVs. The framework suffers from considerable communication overhead and computational burden, and it does not provide flexibility for the agents to have their desired objective function. \citet{LeFloch1} model the scheduling optimization problem as a weakly concave function which is solved iteratively between the agents by the projected gradient method. To have a fully distributed EVCS, the authors use the Lagrangian to consider the coupled constrains in the optimization problem which leads to suboptimal results. In addition, their framework does not consider EVAs' objective function, and there is a strict assumption on DNO's objective function. A hierarchical ADMM for the EV load coordination is proposed by \citet{VERSCHAE2016922} which is derived based on the \textit{sharing problem}. Although their method has mild assumptions on the objective function, it suffers from slow convergence and considerable communication overhead. \citet{Hill1} expand the framework proposed by \citet{Li1} to a generic V2G scheduling framework which includes several layers of EVAs. There is, however, no flexibility in their method to include any other desired scheduling objective, the optimization solved by EVAs to calculate the V2G capacity depends on the number of EVs, and the agents update their optimization variable sequentially. \citet{ZOU2018138} develop a sequential trilayer structure, where the transformers as EVA agents reach a consensus price minimizing the generation cost in the first iterative procedure, and EVs maximize their payoff function according to their transformers' price signal. The convergence of the method, however, is based on strict assumptions on the objective function. 

As discussed above, the desired objective functions of either a lower-level agent, e.g. EV owners, or a higher-level agent, e.g. EVAs, is missing in the bilayer EV charging methods. Therefore, those methods do not satisfy the practical objectives of EVCS. Furthermore, the proposed trilayer EVCS methods in the literature have several drawbacks: (\RNum{1}) they are designed for specific objectives and do not have the flexibility to accommodate any other objective in which case they have to be redesigned (\citet{Ting1}, \citet{Song1}, \citet{Li1}, \citet{Shao1}, \citet{LeFloch1}, \citet{Hill1}); (\RNum{2}) they are not computationally efficient as the agents update their primal variables sequentially at each iteration, e.g. EVAs have to wait for DNO before updating their primal variable (\citet{Ting1}, \citet{Song1}, \citet{Shao1}, \citet{VERSCHAE2016922}, \citet{Hill1}); (\RNum{3}) the convergence and feasibility of the optimization problem solved by the higher-level agents depend on the number of the lower-level agents (\citet{Li1}, \citet{Hill1}); and (\RNum{4}) to solve the EVCS problem in a distributed manner, the coupled constraints relating to feeder capacities are relaxed by the Lagrangian which results in suboptimal solutions (\citet{Song1,LeFloch1}).
     
\paragraph{Contribution of the paper}
In this paper, we expand the research on the trilayer ECI to address their drawbacks by designing a general framework for the interaction among the ECI's agents so that any desired objective for EVCS can be considered.
We propose a hierarchical distributed EVCS, called HDEVCS, meaning that DNO, EVA, and EV agents solve their optimization problem locally while the coupled constraints among the agents are satisfied.
In HDEVCS, DNO communicates only with EVAs, and each EV agent communicates only with its EVA.
To derive HDEVCS, we exploit the configuration of the trilayer ECI and the mathematical properties of the EVCS problem to develop the hierarchical and clustered \textit{exchange problem} which is solved efficiently in a distributed manner by ADMM. The contributions of this paper are summarized as follows:
  

\begin{itemize}
	\item We propose the multi-agent trilayer EV charging framework where each agent has its own objective function which it solves locally. In HDEVCS, the agents do not share their sensitive information with others, so their privacy is preserved.
	Moreover, the optimization problems solved by DNO and EVAs do not depend on the number of system agents. 
	Consequently, the size of their optimization variables does not change by the integration of new agents.
	       
	\item We develop HDEVCS based on the \textit{exchange problem} which is solved efficiently by ADMM. Owing to the \textit{exchange problem}'s properties, the agents have only one primal variable and update it simultaneously at every iteration of HDEVCS, while in the previous methods the agents update their primal variable sequentially. Therefore, HDEVCS does not suffer from significant communication overhead, and it reduces the convergence time and number of iterations considerably. This feature is shown in Section \ref{Simulation Results} by comparing our method with the hierarchical ADMM  (\citet{VERSCHAE2016922}), since it has more flexibility and less strict assumptions compared to the other trilayer methods in the literature.
	
	\item We benefit from the hierarchical configuration of the derived \textit{exchange problem} to satisfy the inequality coupled constraints relating to the EVA feeder capacities in a fully distributed manner. This feature addresses the problem in the previous approaches where the inequality coupled constraints are included in the EVCS problem by the Lagrangian leading to the suboptimal solutions. 
	
	\item As the agents can have any desired objective function, we embed HDEVCS in the receding horizon (RH) feedback control, called RH-HDEVCS, to give the agents flexibility to change their objective function in any RH iteration without notifying other agents. This feature, which is called Plug-and-Play (PnP) and discussed in Subsection \ref{RH-HDEVCS}, addresses the inflexibility of the proposed trilayer methods in the literature which are designed for a specific objective function.   
\end{itemize} 
\paragraph{Synopsis}
The paper is structured as follows: In Section \ref{System Model Description}, the model of EV charging demand and ECI agents are introduced. In Section \ref{Problem Formulation}, a general framework for EVCS is formulated based on the agents' objective functions. In Section \ref{SectionHDEVCS}, the proposed HDEVCS is derived by the clustered \textit{exchange problem}, and its implementation using RH feedback control is discussed. In Section \ref{Simulation Results}, the metrics to evaluate the performance of the proposed framework are defined, numerical simulation results are shown and discussed for two case studies, and our method is compared with the literature to verify its faster convergence. Finally, the paper is concluded in Section \ref{conclusion}. 

\section{System Model Description}  \label{System Model Description}
In this section, after introducing the notation used in the paper, the model of the trilayer ECI's agents, i.e. EV, EVA, and DNO, is presented.
\subsection{Notation} \label{Notations}
We denote the set of natural and real numbers, respectively, by $\mathbb{N}$ and $\mathbb{R}$. The set of positive real values which are less than or equal to one is denoted by $\mathbb{R}^{+}_{\leqslant{1}}$. We use double bracket to show a discrete set, e.g. $\llbracket1,\mathcal{I}\rrbracket=\{1,2,\dots,\mathcal{I}\}$. The vector notation is denoted by 
\begin{align*}
\mathbf{y}(t)=\big(\mathbf{y}_{1}^T(t),\mathbf{y}_{2}^T(t), \dots,\mathbf{y}_{\mathcal{N}}^T(t)\big)^{T},
\end{align*} 
where $\mathbf{y}_{n}(t)=\big(y_{n}^T(t),y_{n}^T(t+1),\dots,y_{n}^T(t+N-1)\big)^{T}, n \in \llbracket1,\mathcal{N}\rrbracket$, is defined componentwise, and $N \in \mathbb{N}$ is a given time horizon. To further simplify the notation, we use $\mathbf{y}:=\mathbf{y}(t)$ and $\mathbf{y}_{n}:=\mathbf{y}_{n}(t)$. The same notion is used for all other variables.

\subsection{EV Charging Infrastructure Model}  \label{EV Charging Model}
\begin{figure}[h]
	\centering \includegraphics[clip,scale=0.48, trim=4.5cm 4.6cm 6.0cm 5.2cm] 
	{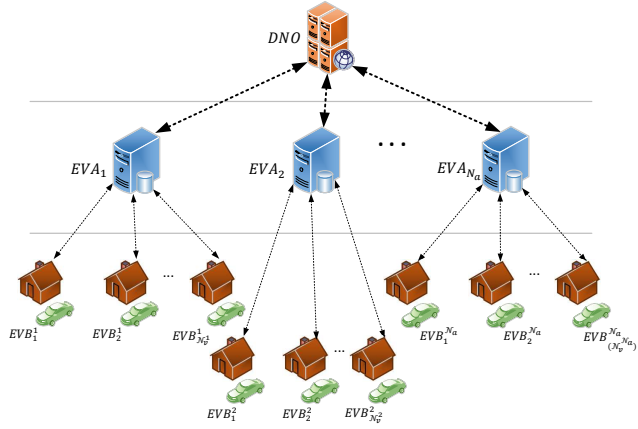}
	\caption{Multi-agent trilayer ECI.}
	\label{fig:testcase}
\end{figure}

In this paper, we consider a trilayer ECI shown in \figref{fig:testcase}.
The system is managed by DNO which aims at optimal EV load coordination in collaboration with EVAs and EVs dispersed throughout the distribution network.
DNO controls the maximum power capacity available for the aggregated charging load demand as well, which is allocated to EVAs for charging EVs at each time instant.
Considering that the charging infrastructure has $\mathcal{N}_{v} \in \mathbb{N}$ EVs, $\mathbb{N}_{a}$ denotes the set of EVAs of the charging infrastructure, $\mathbb{N}_{v}^j$ denotes the set of EVs supplied by $\text{\textit{EVA}}_j$, and their cardinalities are shown by $\mathcal{N}_{a}$ and $\mathcal{N}_{v}^{j}$, respectively. $\text{\textit{EV}}_{i}^j$ shows the $i$th EV supplied by $\text{\textit{EVA}}_j$. 

It is assumed that each $\text{\textit{EV}}_i^j$ is located in either a commercial or a residential building. We use $\text{\textit{EVB}}_{i}^j$ for each $EV_{i}^j$ and its corresponding building, which is modeled as a discrete-time linear system (\citet{8440531,Philipp2}) as
\begin{subequations} \label{RESModel1}
	\begin{align}
	&c_{v_i^j}(t+1)=c_{v_i^j}(t)+T_h p_{v_i^j}(t) \label{RESModel11}\\
	&e_{v_i^j}(t)=d_{v_i^j}(t)+p_{v_i^j}(t) \label{RESModel12},
	\end{align}
\end{subequations}
where $c_{v_i^j}, p_{v_i^j}, e_{v_i^j}, d_{v_i^j} \in \mathbb{R}$, $T_h \in \mathbb{R}^{+}_{\leqslant{1}}$. In \eqref{RESModel11}, $c_{v_i^j}(t)$ is the energy stored in the EV battery at $t$, $d_{v_i^j}(t)$ in \eqref{RESModel12} is the non-EV load demand minus the power generated by the solar panel of the building at $t$, which is also called netload, and $e_{v_i^j}(t)$ is the total load demand (EV charging load plus netload) of the building. $T_h$ in \eqref{RESModel11} is the discretization in time, e.g. $T_h=0.5$ corresponds to 30 min in this paper. $p_{v_i^j}(t)$ is the $\text{\textit{EV}}_i^j$'s charging power and the control variable which is determined by EVCS. Since we assume that each $\text{\textit{EVB}}_i^j$ is provided by a solar panel, and its EV charger has V2G capability, it may supply power to the grid. That is determined by the sign of $e_{v_i^j}(t)$, i.e. $e_{v_i^j}(t)>0$ when $\text{\textit{EVB}}_i^j$ is consuming energy from the grid, and $e_{v_i^j}(t)<0$ when $\text{\textit{EVB}}_i^j$ is feeding energy to the grid.      

The constraints on the EV charging/discharging, related to the charger power rating and the EV battery capacity, are 
\begin{subequations} \label{RESModel2}
	\begin{align}
	&\underline{p}_{v_i^j}(t)\leqslant{p}_{v_i^j}(t)\leqslant{\overline{p}_{v_i^j}}(t) \label{RESModel21}\\
	&\underline{C}_{v_i^j}(t)\leqslant{c_{v_i^j}}(t)\leqslant{\overline{C}_{v_i^j}}(t) \label{RESModel22},
	\end{align}
\end{subequations}
where $\underline{p}_{v_i^j}(t)\in \mathbb{R}$ is the charger's minimum power rating, $\overline{p}_{v_i^j}(t)\in \mathbb{R}$ is the charger's maximum power rating, and $\underline{C}_{v_i^j}(t), \overline{C}_{v_i^j}(t) \in \mathbb{R}$ are the EV battery time-varying constraints which are defined as follows; if $\text{\textit{EV}}_i^j$ is: 
\begin{itemize}
	\item{not plugged in the $\text{\textit{EVB}}_i^j$'s charger, $\underline{C}_{v_i^j}(t)=\overline{C}_{v_i^j}(t)=0$.} 
	\item{plugged in the $\text{\textit{EVB}}_i^j$'s charger, but it is in idle mode, ${\underline{C}_{v_i^j}(t)=0}$ and $\overline{C}_{v_i^j}(t)=C_{v_i^j}$, where $C_{v_i^j} \in \mathbb{R}$ is the maximum EV battery energy capacity.}
	\item{plugged in the $\text{\textit{EVB}}_i^j$'s charger, and it is needed by time $t$, $\underline{C}_{v_i^j}(t)=\overline{C}_{v_i^j}(t)=C_{v_i^j}$, which is according to our assumption that all EVs must be fully charged by their departure time.} 
\end{itemize}

We define the set of feasible charging trajectories of $\text{\textit{EVB}}_i^j$ as
\begin{align}
\label{RESSet}
\mathbb{U}_{v_i^j} =\bigg\{{\mathbf{p}_{v_i^j} \in \mathbb{R}^{N}}\mid{\eqref{RESModel1}-\eqref{RESModel2}~\forall t \in \llbracket k,k+N-1\rrbracket} \bigg\}.
\end{align}

It is the responsibility of each EVA to prevent its aggregated EV power demand from exceeding the available feeder capacity. This is shown by
\begin{subequations} \label{FAConstraint}
	\begin{align}
	&{\mathbf{p}_{a_j}}\leqslant{\overline{\mathbf{P}}_{a_j}},~\overline{\mathbf{P}}_{a_j}=\overline{\mathbf{P}}_{j}-\sum\limits_{i \in \mathbb{N}_{v}^j}^{}{\mathbf{d}_{v_i^j}}  \label{FAConstraint1} \\
	&{\mathbf{p}_{a_j}}=\sum\limits_{i \in \mathbb{N}_{v}^j}^{}{\mathbf{p}_{v_i^j}}\label{FAConstraint2},
	\end{align}
\end{subequations} 
where $\mathbf{p}_{a_j} \in \mathbb{R}^{N}$ is the aggregated EV power demand of $\text{\textit{EVA}}_j$, $\overline{\mathbf{P}}_{a_j}$ is the maximum capacity available for EV power demand, and $\overline{\mathbf{P}}_{j}$ is the maximum capacity of the feeder supplying $\text{\textit{EVA}}_j$. It means that the available capacity to satisfy aggregated EV power demand depends on the total netload demand on the feeder. As we assume that netload demand of the EVBs is not controllable, EVAs can only control the aggregated EV power demand to prevent the total feeder power from exceeding the capacity.

As mentioned before, DNO controls the total EV charging demand so that it does not exceed the maximum available capacity for charging load ($\overline{\mathbf{P}}_{d}$). This constraint is shown by
\begin{subequations}\label{EVMSConstraint}
	\begin{align}
	&-\mathbf{p}_{d}\leqslant{\overline{\mathbf{P}}_{d}} \label{EVMSConstraint1} \\
	&\mathbf{p}_{d}=-\sum\limits_{j \in \mathbb{N}_{a}}^{}{\mathbf{p}_{a_j}} \label{EVMSConstraint2},
	\end{align}
\end{subequations}
where $\mathbf{p}_{d} \in \mathbb{R}^{N}$ is the total EV power demand of the grid. Similar to \eqref{FAConstraint}, DNO can control only the total aggregated EV load demand in grid in order to prevent it from exceeding the capacity. 
\section{Problem Formulation} \label{Problem Formulation}
The ECI shown in \figref{fig:testcase} has a trilayer hierarchical structure including DNO, EVA, and EVB agents. The objective of EVCS is to generate a sequence of feasible EV charging profiles which minimize (maximize) a desirable loss (revenue) function while the operational constraints are satisfied. In this paper, the objective of the trilayer EVCS is threefold which is written as
\begin{equation} \label{ObjFunGeneral1}
\begin{split}
&V:=\min_{\mathbf{p}_{d},\mathbf{p}_{a},\mathbf{p}_{v}}\mathbf{F}_{d}(\mathbf{p}_{d})+\sum\limits_{j \in \mathbb{N}_{a}}^{}{\big(\mathbf{F}_{a_j}(\mathbf{p}_{a_j})+\sum\limits_{i \in \mathbb{N}_{v}^{j}}^{}{\mathbf{F}_{v_i^j}(\mathbf{p}_{v_i^j})}\big)}\\
&=\min_{\mathbf{p}_{d},\mathbf{p}_{a},\mathbf{p}_{v}}\mathbf{F}_{d}(\mathbf{p}_{d})+\mathbf{F}_{a}(\mathbf{p}_{a})+\mathbf{F}_{v}(\mathbf{p}_{v}) \\
&\qquad\qquad \quad \text{s.t.}~\mathbf{p}_{v_i^j} \in \mathbb{U}_{v_i^j},~\forall i \in \mathbb{N}_{v}^{j},~\forall j \in \mathbb{N}_{a}\\ 
&\qquad\qquad\qquad \qquad  \qquad  \eqref{FAConstraint} ~\& ~ \eqref{EVMSConstraint},
\end{split}
\end{equation} 
where $\mathbf{F}_{d}$, $\mathbf{F}_{a}$, and $\mathbf{F}_{v}$ are convex objective functions of DNO, EVAs, and EVBs, respectively. The constraint set of EVBs \eqref{RESSet} as well as \eqref{FAConstraint1} and \eqref{EVMSConstraint1} are the local constraints, while \eqref{FAConstraint2} and \eqref{EVMSConstraint2} are the coupled constraints of the scheduling problem \eqref{ObjFunGeneral1}. 

We define an auxiliary variable for each $\text{\textit{EVA}}_j$ as $\mathbf{p}_{au_j}=-\mathbf{p}_{a_j}$ in order to rewrite \eqref{ObjFunGeneral1} in the \textit{exchange problem} form which can be solved efficiently by ADMM. Considering \eqref{FAConstraint2} and \eqref{EVMSConstraint2}, we have
\begin{equation} \label{Sum0}
\mathbf{p}_{d}+\sum\limits_{j \in \mathbb{N}_{a}}^{}{\mathbf{p}_{a_j}}+\sum\limits_{j \in \mathbb{N}_{a}}^{}{\mathbf{p}_{au_j}}+\sum\limits_{j \in \mathbb{N}_{a}}^{}\sum\limits_{i \in \mathbb{N}_{v}^{j}}^{}{\mathbf{p}_{v_j^j}}=\mathbf{0},
\end{equation}
where $\mathbf{0} \in \mathbb{R}^{N}$ is the zero vector.
This inequality constraint helps us write \eqref{ObjFunGeneral1} as a hierarchical \textit{exchange problem}. In that way, we do not write the local constraints and only show the coupled constraint \eqref{Sum0}.  
\begin{equation} \label{ObjFunGeneral3}
\begin{split}
&V:=\min_{\mathbf{p}_{d},\mathbf{p}_{a},\mathbf{p}_{v}}\mathbf{F}_{d}(\mathbf{p}_{d})+\mathbf{F}_{a}(\mathbf{p}_{a})+\mathbf{F}_{au}(\mathbf{p}_{au})+\mathbf{F}_{v}(\mathbf{p}_{v})\\
& \qquad\qquad\qquad\qquad\qquad\quad \text{s.t.}~ \eqref{Sum0},
\end{split}
\end{equation}
in which $\mathbf{F}_{au}$ denotes the indicator function depending on the EVAs' auxiliary variable ($\mathbf{p}_{au}$), and it is defined as follows for each $\text{\textit{EVA}}_j$
\begin{equation} \label{Fau}
\mathbf{F}_{au_j}(\mathbf{p}_{au_j})=\begin{cases}
0, & \text{if}~\mathbf{p}_{au_j}=-\mathbf{p}_{a_j}\\
\inf, & \text{otherwise}.
\end{cases}
\end{equation}

In the following subsection, we show that the optimization problem in \eqref{ObjFunGeneral3} is the hierarchical \textit{exchange problem} which can be solved in a distributed manner by ADMM. We call it HDEVCS as our proposed framework results in the network through which DNO communicates only with EVAs, and EVAs communicate with their own EVs. That is, there is no direct communication between DNO (the highest-level agent) and EVs (the lowest-level agent), so the proposed EVCS has a hierarchical structure. Nevertheless, all the agents solve their optimization problem locally and simultaneously which decreases the convergence time and communication overhead. 

\section{Hierarchical Distributed Charging Scheduling} \label{SectionHDEVCS}
In this section, the EVCS problem is manipulated mathematically to derive the hierarchical and clustered \textit{exchange problem}. Then, ADMM is applied to solve it in a fully distributed manner.
\subsection{Exchange Problem and HDEVCS} 
To extract HDEVCS framework, we introduce the following aggregated objective function and optimization variable.
\begin{subequations} 
\begin{align*} 
&\mathbf{F} = (\mathbf{F}_{v},\mathbf{F}_{au},\mathbf{F}_{a},\mathbf{F}_{d})\\
&\mathbf{p}= (\mathbf{p}_{v},\mathbf{p}_{au},\mathbf{p}_{a},\mathbf{p}_{d}),
\end{align*}
\end{subequations}
where $\mathbf{F},\mathbf{p} \in \mathbb{R}^{\mathcal{N}_{F}N}$, ${\mathcal{N}_{F}}=({\mathcal{N}_{v}}+2 \times {\mathcal{N}_{a}}+1) \in \mathbb{N}$. We rewrite \eqref{ObjFunGeneral3} as
\begin{equation}  \label{ShP1}
\begin{split}
&V:=\min \sum\limits_{n=1}^{\mathcal{N}_{F}}{\mathbf{F}_{n}(\mathbf{p}_{n})}+\mathbf{G}(\mathbf{z})\\
& \quad \text{s.t.}~\mathbf{p}_{n}=\mathbf{z}_{n},~\forall n \in \llbracket1,\mathcal{N}_{F}\rrbracket,
\end{split}
\end{equation}
where $\mathbf{z} \in \mathbb{R}^{\mathcal{N}_{F}N}$, and $\mathbf{G}(\mathbf{z})$ is the indicator function defined as
\begin{equation} \label{G}
\mathbf{G}(\mathbf{z})=\begin{cases}
0, & \text{if} \sum\limits_{n =1}^{\mathcal{N}_{F}}{\mathbf{z}_{n}}=\mathbf{0}\\
\inf, & \text{otherwise}.
\end{cases}
\end{equation}
The augmented Lagrangian of \eqref{ShP1} is written as
\begin{equation} \label{Lagrangian}
\mathcal{L}_{\rho}(\mathbf{p},\mathbf{z},\pmb{\lambda}) = \mathbf{F}(\mathbf{p})+\mathbf{G}(\mathbf{z})+\pmb{\lambda}^{T}(\mathbf{p} -\mathbf{z})+\frac{\rho}{2}\norm{\mathbf{p} - \mathbf{z}}_{2}^{2},
\end{equation}
in which $\rho$ is the penalty factor, and $\pmb{\lambda}$ is the Lagrangian variable, also known as the dual variable. The Lagrangian \eqref{Lagrangian} can be solved by ADMM which is a variant of augmented Lagrangian approach and uses partial updates of the dual variables in each iteration. The iterative primal and dual updates of ADMM are
\begin{subequations}   
	\begin{align}
	&\mathbf{p}^{k+1}:= \argmin_{\mathbf{p}} \mathcal{L}_{\rho}(\mathbf{p},\mathbf{z}^{k},\pmb{\lambda}^{k}) \label{ADMM_p1}\\
	&\mathbf{z}^{k+1}:= \argmin_{\mathbf{z}} \mathcal{L}_{\rho}(\mathbf{p}^{k+1},\mathbf{z},\pmb{\lambda}^{k}) \label{ADMM_z1}\\
	&\pmb{\lambda}^{k+1}:= \argmax_{\pmb{\lambda}} {\mathcal{L}_{\rho}(\mathbf{p}^{k+1},\mathbf{z}^{k+1},\pmb{\lambda})}, \label{ADMM_y1} 
	\end{align}
\end{subequations}
in which $k$ is the iteration index. Hereafter, we use the scaled form of ADMM (\citet[Chapter~3.1.1]{Boyd1}) where $\pmb{\Lambda}=\pmb{\lambda}/ \rho$. The first step of ADMM \eqref{ADMM_p1} is expanded as
\begin{equation} \label{ADMM_p2}
\mathbf{p}^{k+1}=\sum\limits_{n=1}^{\mathcal{N}_{F}}{\bigg(\mathbf{F}_n(\mathbf{p}_n)+\frac{\rho}{2}\norm{\mathbf{p}_n-\mathbf{z}_n^{k}+\pmb{\Lambda}_n^{k}}\bigg)},
\end{equation}
which is separable and can be solved in parallel by each agent of ECI, i.e. DNO, EVAs and EVBs. Further details will be provided later.

According to \eqref{ShP1}, each $\mathbf{z}_{n}$ is equivalent to $\mathbf{p}_{n}$. Using the defined auxiliary variable $\mathbf{p}_{au}$ and also \eqref{FAConstraint2} and \eqref{EVMSConstraint2}, we can partition the charging infrastructure into $\mathcal{N}_{c}=(\mathcal{N}_{a}+1)$ clusters denoted by $\text{\textbf{CL}}$s. For each cluster, the coupled equality constraint which is the summation of the involved agents' primal variables is defined. That equality constraint is equal to zero, and it is known as the \textit{equilibrium constraint}. The clusters and the corresponding \textit{equilibrium constraints} of HDEVCS are defined as follows:
\begin{itemize}
	\item $\text{\textbf{CL}}_{j}$, $\forall j \in \llbracket1,(\mathcal{N}_{c}-1)\rrbracket$, includes $\text{\textit{EVA}}_j$ and $\text{\textit{EVB}}_i^j$, $\forall i \in \mathbb{N}_v^j$, and its \textit{equilibrium constraint} is $\sum\limits_{i=1}^{\mathcal{N}_{v}^j}{\mathbf{p}_{v_i^j}}+\mathbf{p}_{{au}_j}=\mathbf{0}$.
	\item $\text{\textbf{CL}}_{\mathcal{N}_c}$ includes  DNO and all EVAs, and its \textit{equilibrium constraint} is $\sum\limits_{j=1}^{\mathcal{N}_{a}}{\mathbf{p}_{a_j}}+\mathbf{p}_{d}=\mathbf{0}$.
	\item According to the equality constraints of the defined clusters, we have:
\begin{equation}  \label{EVAsets} 
\sum\limits_{n \in \text{\textbf{CL}}_{j}}^{}{\mathbf{z}_n}=\mathbf{0},~\forall j \in \llbracket1,\mathcal{N}_{c}\rrbracket. 
\end{equation}  
\end{itemize}
The above three statements mean that the charging infrastructure includes $\mathcal{N}_c$ clusters where the interaction among the agents within each cluster can be written as the \textit{exchange problem}. To show that, we write the Lagrangian for $\mathbf{z}$ \eqref{ADMM_z1} in a partitioned form as 
\begin{equation}  \label{ADMM_z2} 
\mathcal{L}(\mathbf{z},\pmb{\upsilon}) = \sum\limits_{j=1}^{\mathcal{N}_{c}}\sum\limits_{n \in \text{\textbf{CL}}_{j}}^{}\bigg(\frac{\rho}{2}\norm{\mathbf{p}_n^{k+1}-\mathbf{z}_n+\pmb{\Lambda}_n^k}_2^2+\pmb{\upsilon}_j\mathbf{z}_n\bigg),
\end{equation}
where $\pmb{\upsilon}_j$ is the Lagrangian multiplier corresponding to the coupled equality constraint of $\text{\textbf{CL}}_{j}$. Using the KKT conditions for each cluster and adding up the gradients of Lagrangian in terms of $\mathbf{z}_{n},~n\in \text{\textbf{CL}}_{j}$, we have 
\begin{equation} \label{proof1}
\nabla_{\mathbf{z}_n}\mathcal{L}(\mathbf{z}_n,\pmb{\upsilon}_j)=0\Rightarrow \mathbf{z}_n=\mathbf{p}_n^{k+1}+\pmb{\Lambda}_n^{k}-\frac{1}{\rho}\pmb{\upsilon}_j
\end{equation}
\begin{equation} \label{proof12}
\nabla_{\upsilon_j}\mathcal{L}(\mathbf{z}_n,\pmb{\upsilon}_j)=0\Rightarrow \sum\limits_{n\in \text{\textbf{CL}}_{j}}^{}{\mathbf{z}_n}=\mathbf{0}
\end{equation}
\begin{equation*}
\eqref{proof1},\eqref{proof12}\Rightarrow \sum\limits_{n\in \text{\textbf{CL}}_{j}}^{}{\big(\mathbf{p}_n^{k+1}+\pmb{\Lambda}_n^{k}\big)}-\frac{\mathcal{N}_{c_j}.\pmb{\upsilon}_j}{\rho}=\mathbf{0},
\end{equation*}
where $\mathcal{N}_{c_j}$ denotes the cardinality of $\text{\textbf{CL}}_{j}$. Finally, the dual variable of $\text{\textbf{CL}}_{j}$ is obtained by
\begin{equation} \label{proof3}
\pmb{\upsilon}_j = \rho\big(\overline{\mathbf{p}}_j^{k+1}+\overline{\pmb{\Lambda}}_j^{k}\big),
\end{equation}
in which $\overline{\mathbf{p}}_j=\frac{1}{\mathcal{N}_{c_j}}\sum\limits_{n\in \text{\textbf{CL}}_{j}}^{}{\mathbf{p}_n}$, and $\overline{\pmb{\Lambda}}_j=\frac{1}{\mathcal{N}_{c_j}}\sum\limits_{n\in \text{\textbf{CL}}_{j}}^{}{\pmb{\Lambda}_n}$. Using \eqref{proof1} and \eqref{proof3}, we can find a closed form to update $\mathbf{z}_n \in \text{\textbf{CL}}_{j}$ as
\begin{equation}  \label{ADMM_z3} 
\mathbf{z}_{n}^{k+1}=\mathbf{p}_n^{k+1}-\overline{\mathbf{p}}_j^{k+1}+\pmb{\Lambda}_n^{k}-\overline{\pmb{\Lambda}}_j^{k}.
\end{equation}
Using the gradient method, $\pmb{\Lambda}_n$ is obtained by
\begin{equation} \label{ADMM_y2}
\pmb{\Lambda}_n^{k+1}=\pmb{\Lambda}_n^{k}+\mathbf{p}_n^{k+1}-\mathbf{z}_n^{k+1}.
\end{equation}
Substituting \eqref{ADMM_z3} for \eqref{ADMM_y2} gives
\begin{equation} \label{ADMM_y3}
\pmb{\Lambda}_n^{k+1}=\overline{\mathbf{p}}_j^{k+1}+\overline{\pmb{\Lambda}}_j^{k},
\end{equation}
meaning that the dual updates for all the agents in $\text{\textbf{CL}}_{j}$ shown by $\overline{\pmb{\Lambda}}_j$ are equal and independent of the number of agents. Therefore, $\mathbf{z}_{n}^{k+1}$ is calculated by
\begin{equation}  \label{ADMM_z4} 
\mathbf{z}_{n}^{k+1}=\mathbf{p}_n^{k+1}-\overline{\mathbf{p}}_j^{k+1}.
\end{equation}

In \eqref{ADMM_p2}, substituting \eqref{ADMM_z4} for $\mathbf{z}_n$ eliminates the second primal-update step of ADMM \eqref{ADMM_z1}, so there is not any sequential primal update in HDEVCS owing to reformulating EVCS as the \textit{exchange problem}. Now, we can show how primal variables ($\mathbf{p}_n$) are updated in parallel by EVBs, EVAs, and DNO at each iteration of HDEVCS using ADMM.

-- primal variable update for $\textit{EVB}_i^j$, $\forall i \in \text{\textbf{CL}}_{j}$:
\begin{equation}\label{EVPUpdate}
\begin{split}
&\mathbf{p}_{v_i^j}^{k+1} = \argmin_{\mathbf{p}_{v_i^j}} \big(\mathbf{F}_{v_i^j}(\mathbf{p}_{v_i^j})+\frac{\rho}{2}\norm{\mathbf{p}_{v_i^j}-\mathbf{p}_{v_i^j}^{k}+\overline{\mathbf{p}}_{j}^{k}+\overline{\pmb{\Lambda}}^{k}_{j}}_2^2 \big)\\
& \qquad\qquad\qquad\qquad\qquad\quad \text{s.t.}~\eqref{RESSet}.
\end{split}
\end{equation}

-- primal variable update for $\textit{EVA}_j$, $j \in \text{\textbf{CL}}_{j}$:
\begin{equation}\label{EAPUpdate}
\begin{split}
&\mathbf{p}_{a_j}^{k+1} = \argmin_{\mathbf{p}_{a_j},\mathbf{p}_{au_j}} \big(\mathbf{F}_{a_j}(\mathbf{p}_{a_j})+\frac{\rho}{2}\norm{-\mathbf{p}_{a_j}+\mathbf{p}_{a_j}^{k}+\overline{\mathbf{p}}_{j}^{k}+\overline{\pmb{\Lambda}}^{k}_{j}}_2^2\\
&+\frac{\rho}{2}\norm{\mathbf{p}_{a_j}-\mathbf{p}_{a_j}^{k}+\overline{\mathbf{p}}_{\mathcal{N}_{c}}+\overline{\pmb{\Lambda}}^{k}_{\mathcal{N}_{c}}}_2^2 \big)\\
& \quad\qquad\qquad\qquad\qquad\qquad\quad \text{s.t.} ~\eqref{FAConstraint},
\end{split}
\end{equation}
in which $-\mathbf{p}_{a_j}$ substitutes for $\mathbf{p}_{au_j}$ in the second right-hand-side expression, meaning $\text{\textit{EVA}}_j$ deals with only one primal variable.  

-- primal variable update for DNO:
\begin{equation}\label{DNOPUpdate}
\begin{split}
&\mathbf{p}_{d}^{k+1} = \argmin_{\mathbf{p}_{d}} \big(\mathbf{F}_{d}(\mathbf{p}_{d})
+\frac{\rho}{2}\norm{\mathbf{p}_{d}-\mathbf{p}_{d}^{k}+\overline{\mathbf{p}}_{\mathcal{N}_c}^{k}+\overline{\pmb{\Lambda}}^{k}_{\mathcal{N}_c}}_2^2 \big)\\
& \qquad\qquad\qquad\qquad\qquad\quad \text{s.t.}~\eqref{EVMSConstraint}.
\end{split}
\end{equation}

After updating the primal variables by all the agents in parallel \eqref{EVPUpdate}-\eqref{DNOPUpdate}, EVBs send out their updated variable to their EVAs, and EVAs transmit their updated variable to DNO. The average power ($\overline{\mathbf{p}}_j$) and the dual variable ($\overline{\pmb{\Lambda}}_j$) for each $\textbf{CL}_j$ are updated by EVAs and DNO. To lower the communication overhead, $\overline{\pmb{\Omega}}_j$ is broadcast to the other agents in $\text{\textbf{CL}}_{j}$ which is defined as $\overline{\pmb{\Omega}}_j=\overline{\mathbf{p}}_j+\overline{\pmb{\Lambda}}_j$. The whole procedure of the proposed HDEVCS is shown in Algorithm \ref{algorithm1}, where $r$ and $s$ are primal and dual residuals, respectively, and $th_p$ and $th_d$ are their corresponding feasibility tolerance. For more details about primal and dual residuals and stopping criteria, we refer to Appendix and \citet[Chapter~3.3]{Boyd1}.%
\begin{algorithm} 
	\caption{HDEVCS}
	\begin{algorithmic}[1]
		\WHILE{$err_p>th_p$ or $err_d>th_d$}
		\ForPA{EVBs, EVAs and DNO}
		\STATE Update $\mathbf{p}_{v_i^j}^{k+1}~\forall \text{\textit{EVB}}_i^j,~i \in \mathbb{N}_v^j,~j \in \mathbb{N}_a,$ by \eqref{EVPUpdate}.
		\STATE Update $\mathbf{p}_{a_j}^{k+1} ~\forall \text{\textit{EVA}}_j,~j \in \mathbb{N}_a$ by \eqref{EAPUpdate}.
		\STATE Update $\mathbf{p}_{d}^{k+1}$ by \eqref{DNOPUpdate}.
		\EndPA
		\ForP{$j=1:\mathcal{N}_{c}$}
		\IF{$j \in \llbracket1,\mathcal{N}_{c}-1\rrbracket$}
		\STATE $\text{\textit{EVA}}_j$ receives $\mathbf{p}_{v_i^j}^{k+1}, \forall$ $\text{\textit{EVB}}_i^j \in \text{\textbf{CL}}_{j}$. 
		\STATE Update $\overline{\mathbf{p}}_j^{k+1}=\frac{1}{\mathcal{N}_{c_j}}\sum\limits_{i\in \mathbb{N}_v^j}^{}{\mathbf{p}_{v_i^j}^{k+1}}$.
		\STATE Update $\overline{\pmb{\Lambda}}_j^{k+1}=\overline{\mathbf{p}}_j^{k+1}+\overline{\pmb{\Lambda}}_j^{k}$.
		\STATE Broadcast $\overline{\pmb{\Omega}}_j^{k+1}$ to $\forall \text{\textit{EVB}}_i^j \in \text{\textbf{CL}}_{j}$.
		\STATE Update $\mathbf{r}_j^{k+1}$, $\mathbf{s}_{au_j}^{k+1}$ and $\mathbf{s}_{v_i^j}^{k+1}$, $\forall i \in \mathbb{N}_v^j$.
		\ELSE 
		\STATE DNO receives $\mathbf{p}_{a_j}^{k+1}, \forall j \in \mathbb{N}_{a}$. 
		\STATE Update $\overline{\mathbf{p}}_{\mathcal{N}_{c}}^{k+1}=\frac{1}{(\mathcal{N}_{a}+1)}(\sum\limits_{j \in \mathbb{N}_a}^{}{\mathbf{p}_{a_j}^{k+1}}+\mathbf{p}_{d}^{k+1})$.
		\STATE Update $\overline{\pmb{\Lambda}}_{\mathcal{N}_{c}}^{k+1}=\overline{\mathbf{p}}_{\mathcal{N}_{c}}^{k+1}+\overline{\pmb{\Lambda}}^{k}_{\mathcal{N}_{c}}$.
		\STATE Broadcast $\overline{\pmb{\Omega}}_{\mathcal{N}_c}^{k+1}$ to all EVAs. 
		\STATE Update $\mathbf{r}_{j}^{k+1}$, $\mathbf{s}_d^{k+1}$ and , $\mathbf{s}_{a_i}^{k+1}$, $\forall i \in \mathbb{N}_a$.
		\ENDIF
		\EndP
		\ENDWHILE
	\end{algorithmic} 
	\label{algorithm1}
\end{algorithm}

\begin{figure}[h]
	\centering \includegraphics[clip,scale=0.48, trim=4.5cm 5.2cm 6.0cm 5.2cm] 
	{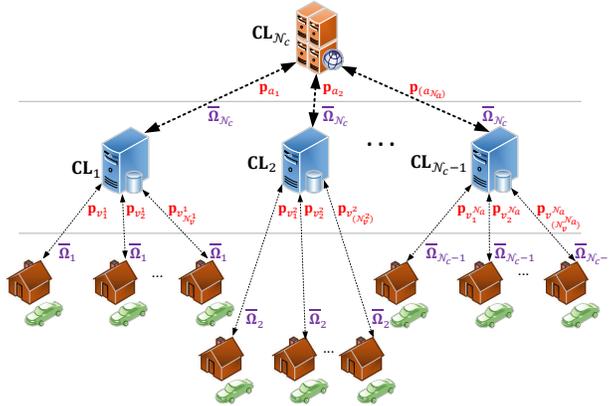}
	\caption{HDEVCS's communication network and broadcast signals.}
	\label{fig:testcase2}
\end{figure}

The communication network links between the agents as well as the broadcast variables within each $\text{\textbf{CL}}_{j}$ are shown in \figref{fig:testcase2}. Note that the index of each cluster, i.e. $\text{\textbf{CL}}_{j}$, is shown next to the agent which updates the average power ($\overline{\mathbf{p}}_j$) and the dual variable ($\overline{\pmb{\Lambda}}_j$).  
\subsection{Receding Horizon HDEVCS: RH-HDEVCS} \label{RH-HDEVCS}
In this subsection, we apply receding horizon feedback control to HDEVCS to generate the optimal control sequences for the EV chargers.
Utilizing RH gives significant flexibility to the EV owners, EVAs, and DNO to change their objective functions, as long as it is convex and feasible, at any RH iteration. In that case, as each agent's optimization function is solved locally, they do not need to notify other agents about changing the objective function. This feature is called PnP (\citet{Philipp2}).

RH-HDEVCS shown in Algorithm \ref{algorithm2} executes the following steps consecutively. First, according to the requested EV charging and total netload demands of the system at time instant $t$ over the time horizon $N \in \mathbb{N}$, Algorithm \ref{algorithm1} finds the optimal values of the primal variable indicated by $\mathbf{p}^\star(\cdot)$. Then, only the first element of $\mathbf{p}^\star(\cdot)$ is implemented by each agent. Lastly, $t$ is incremented, and the same procedure is repeated.
\begin{algorithm}
	\caption{RH-HDEVCS}
	\begin{algorithmic}[1]
		\WHILE{$t\leqslant N$}
		\STATE Update EVs' arrival/departure time, initial state of energy, and charging energy demand, and EVAs' and DNO's available capacities.
		\STATE Run HDEVCS, Algorithm \ref{algorithm1}, to calculate $\mathbf{p}^\star(\cdot)$.
		\STATE Apply the first element of the optimal value of $\mathbf{p}^\star(\cdot)$ for each agent.
		\STATE Increment the time index $t$. 
		\ENDWHILE
	\end{algorithmic}
	\label{algorithm2}
\end{algorithm}
\section{Numerical Simulation} \label{Simulation Results} 
In this section, the performance of HDEVCS is evaluated for two case studies, a small-scale system and a large-scale system, which are called \textit{System\textsubscript{1}} and \textit{System\textsubscript{2}}, respectively. To show the effectiveness of the proposed EVCS, its performance is compared with uncoordinated (\textit{uCC}) and semi-coordinated (\textit{sCC}) charging methods. In \textit{uCC}, EVs start charging with the maximum power rating ($\overline{\mathbf{p}}_{v}$) as soon as they are plugged in, while in \textit{sCC} they charge with a constant power rating while they are plugged in. 

Through all the simulations, the maximum power rating for EVBs is $4$ kW. The netload dataset of EVBs is collected from the Australian electricity company-Ausgrid \cite{Ausgrid} provided for $300$ residential customers, and the wholesale price is available from the California Independent System Operator-CAISO \cite{CAISO}. All the simulations are executed by MATLAB on a PC with Intel\textsuperscript{$\text{\textregistered}$} Core\textsuperscript{\texttrademark} i$7-4770$ $3.40$ GHz CPU, $4$ cores and $8$ GB RAM, and the convex optimization problems are solved by CVX \cite{CVX}.   
\paragraph{Evaluation Metrics} 
We define three metrics to assess the performance of HDEVCS for LVM. The first metric is the peak-to-peak (PTP) value of the aggregated power demand (neatload plus EV load) seen by DNO, which is the difference between the maximum and minimum power demand over the time horizon $N \in \mathbb{N}$. PTP is defined as
\begin{equation} \label{PTP} 
\text{PTP} = \max_{t \in \llbracket1,N\rrbracket} \varepsilon (t) - \min_{t \in \llbracket1,N\rrbracket} \varepsilon (t),
\end{equation}
where $\varepsilon(t) = \sum\limits_{j \in \mathbb{N}_{a}}^{} \sum\limits_{i \in \mathbb{N}_v^j}^{}e_{v_i^j}(t)$, $\varepsilon(t) \in \mathbb{R}$. The second performance metric is the peak-to-average (PTA) value (\citet{Monti1}), which is given by
\begin{equation} \label{PTA} 
\text{PTA} = \frac{\max_{t \in \llbracket1,N\rrbracket} \varepsilon (t)}{\overline{E}(t)},
\end{equation}
where $\overline{E}(t) \in \mathbb{R}$ is the average of the aggregated netload demand of grid over the time horizon $N \in \mathbb{N}$, and it is calculated as $\overline{E}(t)= \frac{1}{N}\sum\limits_{t \in \llbracket1,N\rrbracket}^{}\sum\limits_{j \in \mathbb{N}_{a}}^{} \sum\limits_{i \in \mathbb{N}_v^j}^{}d_{v_i^j}(t)$. The last metric is the root-mean-square (RMS) deviation from the average power demand which is given by
\begin{equation} \label{RMS} 
\text{RMS} = \sqrt{\frac{1}{N}\sum\limits_{t=1}^{N}{\big(\overline{E}(t)-\varepsilon(t)\big)^2}}. 
\end{equation}

\begin{figure}[t]
	\centering
	\subfloat 
	{\includegraphics[clip,scale=0.53, trim=2.3cm 9.2cm 9.5cm 11.0cm]{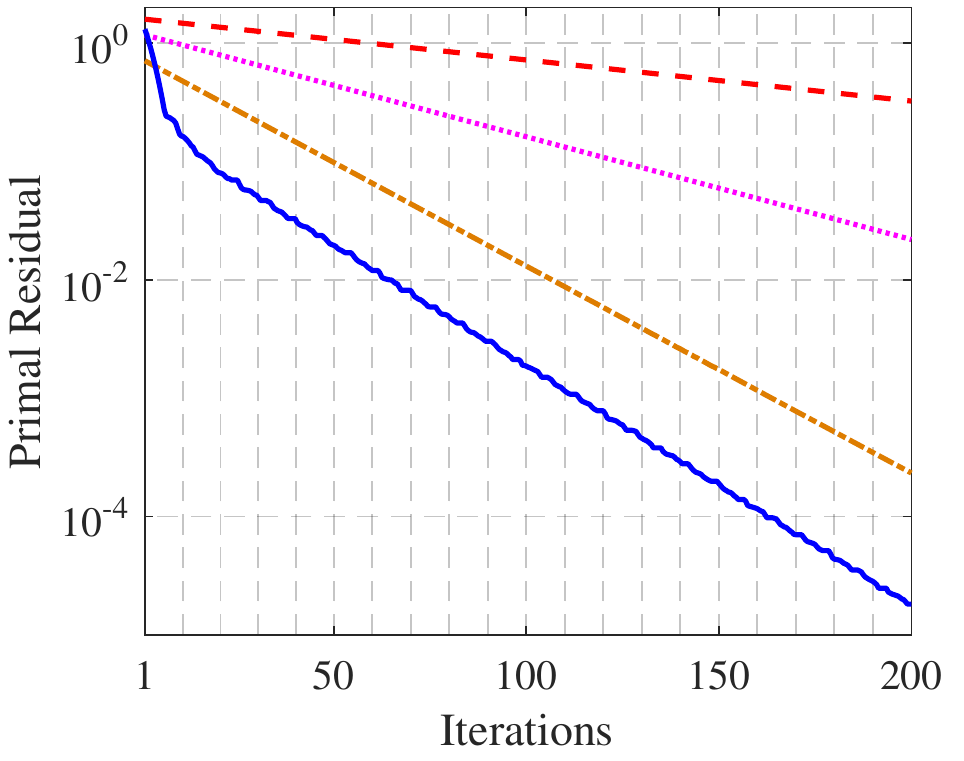}\label{fig:Primal}}
	\subfloat 
	{\includegraphics[clip,scale=0.53, trim=2.3cm 9.2cm 7.7cm 11.0cm]{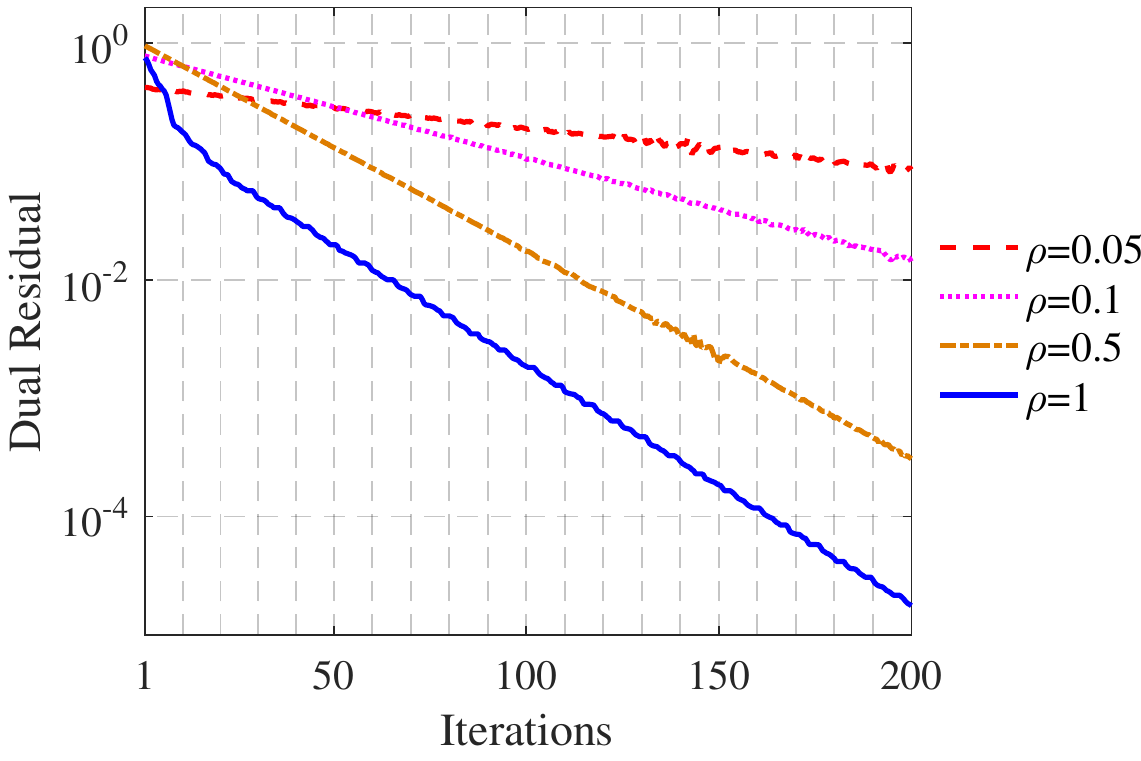}\label{fig:Dual}}
	\caption{Convergence of primal (left) and dual (right) residuals for different $\rho$ values.}\label{fig:rhoanalysis}
\end{figure}

\paragraph{Convergence Analysis}
To investigate the convergence behavior of ADMM for HDEVCS, we run the numerical simulation using four different values of the penalty factor ($\rho$).
As shown in \figref{fig:rhoanalysis}, both the primal and the dual residuals linearly converge to zero for smaller values of $\rho$.
For the largest penalty factor value (i.e. $\rho=1$ in our case), after $20$ iterations the residuals linearly decrease, although their initial rate of convergence is faster.
Also after $200$ iterations, the residuals are still considerable for $\rho \in \{0.05,0.1\}$ which confirms their slow rate of convergence. In addition, while the primal residual chatters for only $\rho=1$, the dual residual chatters for all simulated values of the penalty factor. Nevertheless, both residuals' convergence is acceptable, and they reach a desired value within around $50$ iterations. According to the results, we pick  $\rho=1$ to run the numerical simulations in the next subsections. It is worthwhile to mention that the convergence rate can be improved by the adaptive penalty term which is not used in this paper.

\subsection{\textit{System\textsubscript{1}}: Small-Scale Case Study}
\textit{System\textsubscript{1}} consists of $5$ EVAs each of which is supplying $60$ EVs. The simulations are executed for two scenarios with different set of objective functions for EVAs and EVs. Based on National Household Travel Survey (NHTS) 2017 (\citet{NHTS}), the required charging energy, arrival time, and departure time of the EVs are generated as follows: the initial and designated EVs' battery energies are normally distributed over $[8,10]$ kWh and $[22,25]$ kWh, respectively; EVs' arrival and departure times are normally distributed in $[\texttt{16:30},\texttt{20:30}]$ and $[\texttt{6:00},\texttt{9:30}]$, respectively.  
\paragraph{Scenario 1}
In the first scenario, the purpose of DNO is to minimize the peak load demand in the system which is equivalent to LVM (\citet{Sortomme1}). Defining the aggregated netload demand of $\textit{EVA}_j$, $j \in \mathbb{N}_a$ by \eqref{FAAvg}
\begin{align}
\label{FAAvg}
&d_{a_j}(t):= \sum\limits_{i \in \mathbb{N}_v^j}^{} d_{v_i^j}(t),
\end{align} 
LVM is obtained by minimizing
\begin{align} \label{DNOObjS1}
\begin{split}
&\mathbf{F}_d(\mathbf{p}_d) := \big(\overline{\textbf{E}}-\mathbf{p}_{d}-\sum\limits_{j \in \mathbb{N}_a}^{} \textbf{d}_{a_j}\big)^{2} \\
&   \quad \qquad \qquad 
\text{s.t.}~\eqref{EVMSConstraint}.
\end{split}
\end{align}
While we assume that EVAs' purpose is only to keep their aggregated EV charging demand less than the feeder capacity constraint \eqref{FAConstraint1}, EVs aim at reducing their charging cost. Accordingly, $\mathbf{F}_{a_j}$ is defined by
\begin{align} \label{EVAObjS1}
\begin{split}
&\mathbf{F}_{a_j}(\mathbf{p}_{a_j}) := \mathbf{I}_{a_j}(\mathbf{p}_{a_j})\\
&   \quad \text{s.t.}~\eqref{FAConstraint},~\forall j \in \mathbb{N}_a,
\end{split}
\end{align}
where $\mathbf{I}_{a_j}$ is the indicator function which is defined as
\begin{equation} \label{I}
\mathbf{I}_{a_j}(\mathbf{p}_{a_j})=\begin{cases}
0, & \text{if}~\eqref{FAConstraint}~\text{is satisfied}\\
\inf, & \text{otherwise}.
\end{cases}
\end{equation}
$\text{\textit{EV}}_{v_i^j}$'s objective function is defined by
\begin{align} \label{EVBObjS1}
\begin{split}
&\mathbf{F}_{v_i^j}(\mathbf{p}_{v_i^j}):=\mathbf{\Pi}^{T}.\mathbf{p}_{v_i^j}\\
&\hspace{-8pt}\text{s.t.}~\eqref{RESSet},~\forall i \in \mathbb{N}_{v}^j,~\forall j \in \mathbb{N}_a,
\end{split}
\end{align} 
where $\mathbf{\Pi} \in \mathbb{R}^N$ is the electricity wholesale price. 

\begin{figure}[h!]
	\centering \includegraphics[clip,scale=0.45, trim=1cm 9.2cm 0cm 10.0cm] 
	{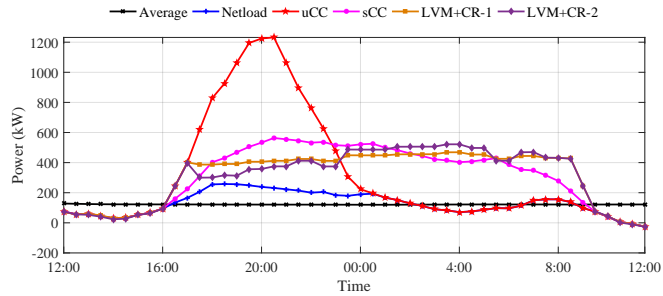}
	\caption{Aggregated load profile for \textit{uCC}, \textit{sCC} and HDEVCS for two different CR weighting factors, $1$ and $10$.}
	\label{fig:uCC_CCMPC_CCnMPC}
\end{figure}

The simulation results of HDEVCS compared with \textit{uCC} and \textit{sCC} are illustrated in \figref{fig:uCC_CCMPC_CCnMPC}. \textit{LVM+CR-1} and \textit{LVM+CR-2} show the aggregated load profile (the aggregated EV charging demand plus netload) when the weighting factor for \eqref{EVBObjS1} is $1$ and $10$, respectively. By decreasing the weighting factor, the load profile becomes smoother at the cost of a more expensive charging for the EV owners.

\begin{figure}[h!]
	\centering
	\subfloat 
	{\includegraphics[clip,scale=0.53, trim=2.3cm 9.2cm 9.3cm 11.0cm]{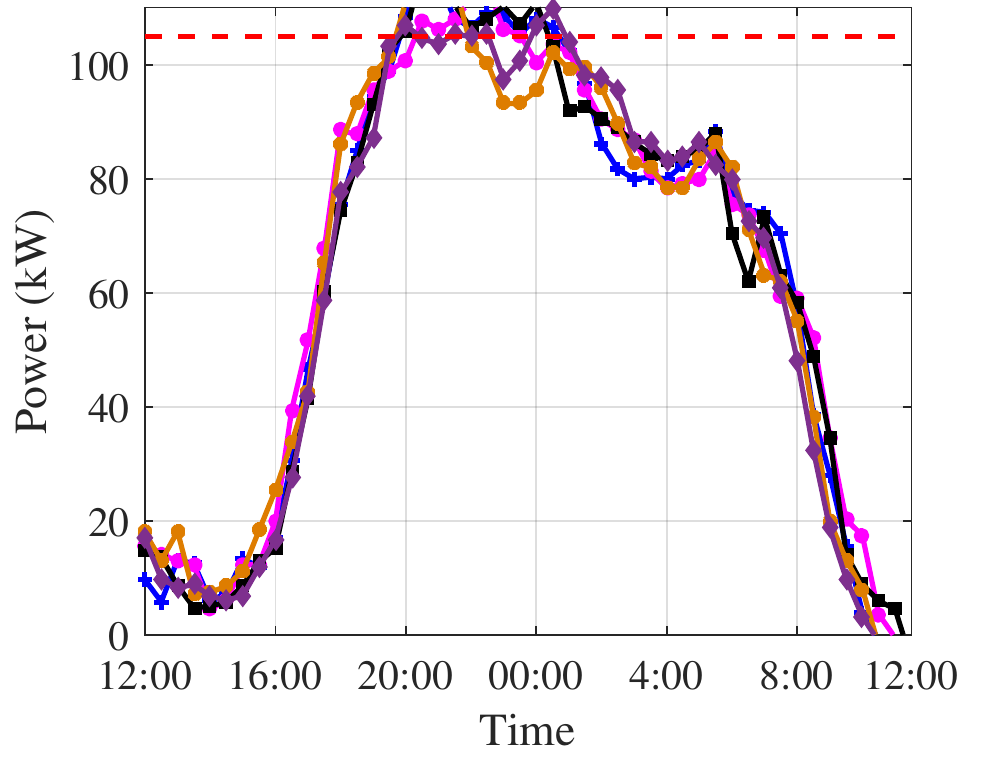}\label{CVio}}
	\subfloat 
	{\includegraphics[clip,scale=0.53, trim=1.5cm 9.2cm 8.0cm 11.0cm]{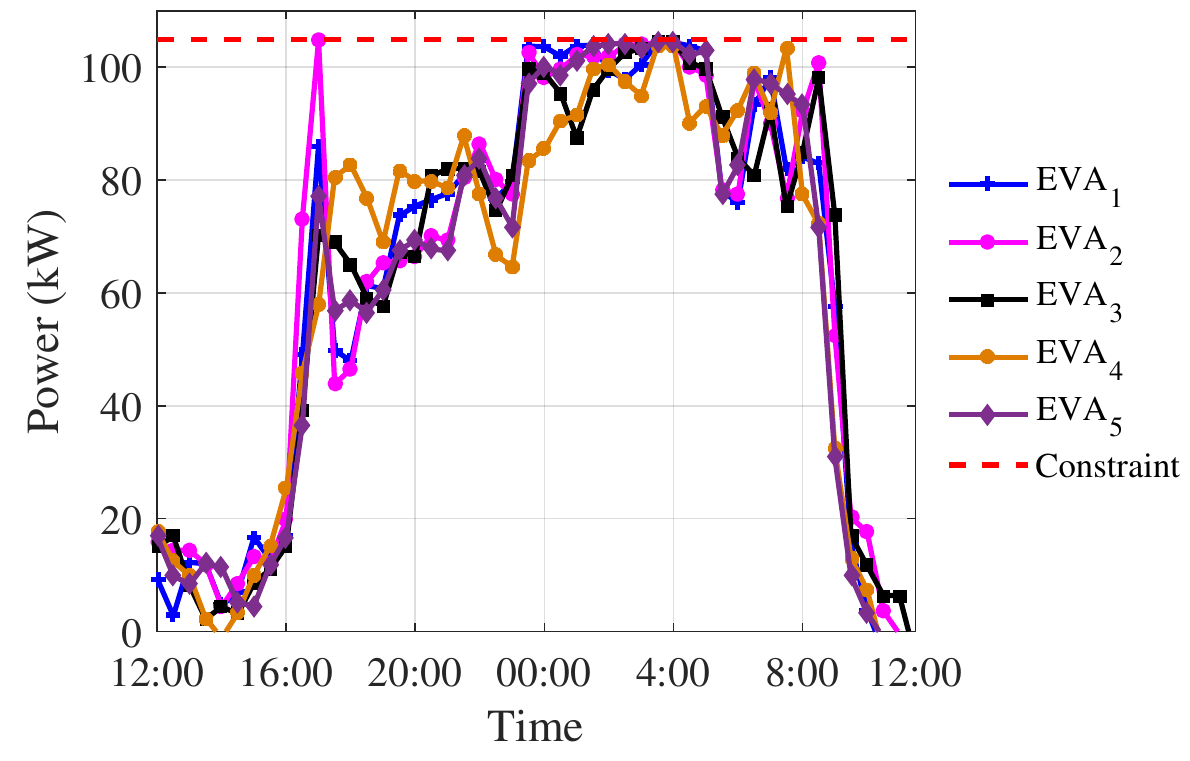}\label{CSat}}
	\caption{EVAs' feeder constraint violation by \textit{sCC} (left) and feeder constraint satisfaction by HDEVCS (right) for \textit{LVM+CR-2} mode.}\label{CCompar}
\end{figure}

In \figref{CCompar}, it is shown that EVCS is effective in limiting the aggregated load demand to the capacity of the EVAs' feeders. While \textit{sCC} does not have any control on the aggregated EV load (\figref{CCompar}: left), HDEVCS does not let the capacity constraints ($105$ kW) be violated (\figref{CCompar}: right) even if the EV agents choose a high weighting factor to greedily reduce their charging cost. These results highlight the importance of optimal coordination of EV charging in supplying more load demand without expansion of the grid capacity. 
\paragraph{Scenario 2}
In the second scenario, we assume that DNO's purpose is still LVM \eqref{DNOObjS1}, while EVAs aim at reducing the aggregated charging cost, and EVs plan for reducing their battery degradation. Thus, we have
\begin{align} \label{EVAObjS2}
\begin{split}
&\mathbf{F}_{a_j}(\mathbf{p}_{a_j}):=\mathbf{\Pi}^{T}.\mathbf{p}_{a_j}\\
& \quad \text{s.t.}~\eqref{FAConstraint},~\forall j \in \mathbb{N}_a.
\end{split}
\end{align}   
To define EVs' optimization function, we borrow BDR model proposed by \citet{Ma2}
\begin{align} \label{EVBObjS2}
\begin{split}
&\mathbf{F}_{v_i^j}(\mathbf{p}_{v_i^j}):=\gamma_{v_i^j}^1\mathbf{p}_{v_i^j}^2+\gamma_{v_i^j}^2\mathbf{p}_{v_i^j}+\gamma_{v_i^j}^3\\
&\hspace{16pt}\text{s.t.}~\eqref{RESSet},~~\forall i \in \mathbb{N}_{v}^j,~\forall j \in \mathbb{N}_a,
\end{split}
\end{align} 
where $\gamma_{v_i^j}^1$, $\gamma_{v_i^j}^2$, and $\gamma_{v_i^j}^3$ are the constant coefficients depending on the number, nominal voltage value and price of energy units of the battery cells. 

\begin{figure}[h!]
	\centering \includegraphics[clip,scale=0.45, trim=1cm 9.2cm 0cm 10.0cm] 
	{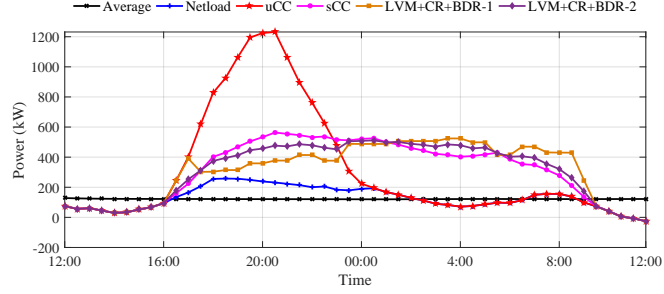}
	\caption{Aggregated load profile for \textit{uCC}, \textit{sCC} and HDEVCS for two different BDR weighting factors, $1$ and $10$.}
	\label{fig:Scenario2LP1}
\end{figure}

The aggregated load profiles obtained by HDEVCS are compared with \textit{uCC} and \textit{sCC} in \figref{fig:Scenario2LP1}. Similar to the first scenario, we run the simulations for two different weighting factors of BDR, i.e. $1$ and $10$, which are shown by \textit{LVM+CR+BDR-1} and \textit{LVM+CR+BDR-2}, respectively. In both modes, the weighting factor of CR which is EVAs' objective function is equal to $10$.

\begin{figure}[h!]
	\centering
	\subfloat 
	{\includegraphics[clip,scale=0.53, trim=2.3cm 9.2cm 9.3cm 11.0cm]{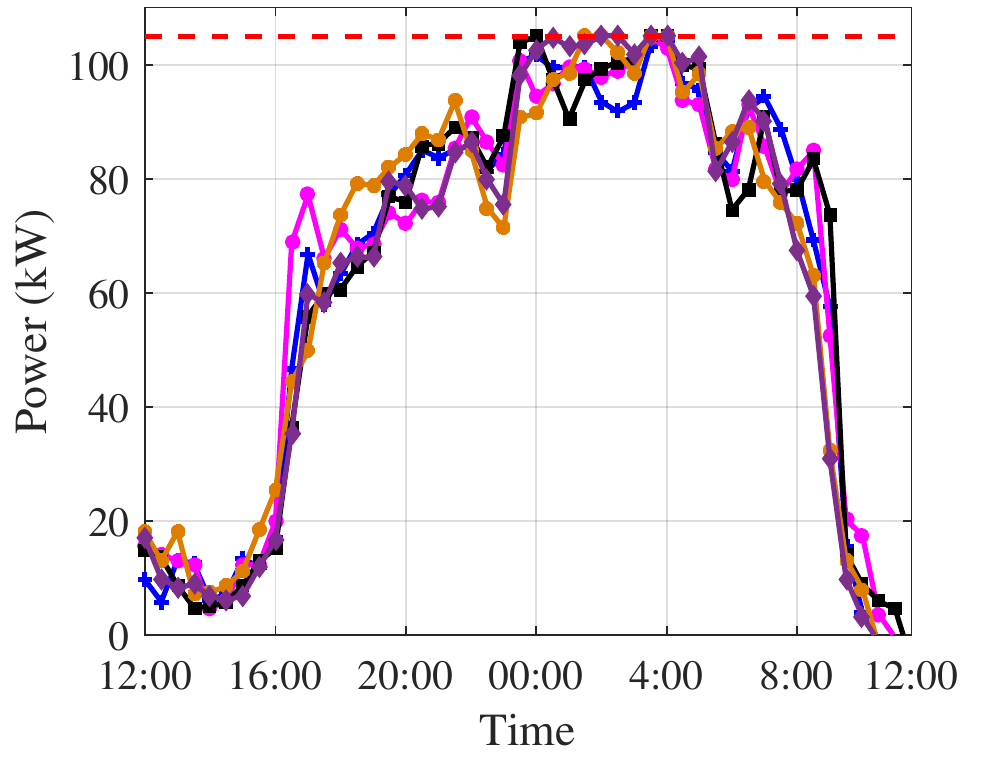}\label{CVio2}}
	\subfloat 
	{\includegraphics[clip,scale=0.53, trim=1.5cm 9.2cm 8.0cm 11.0cm]{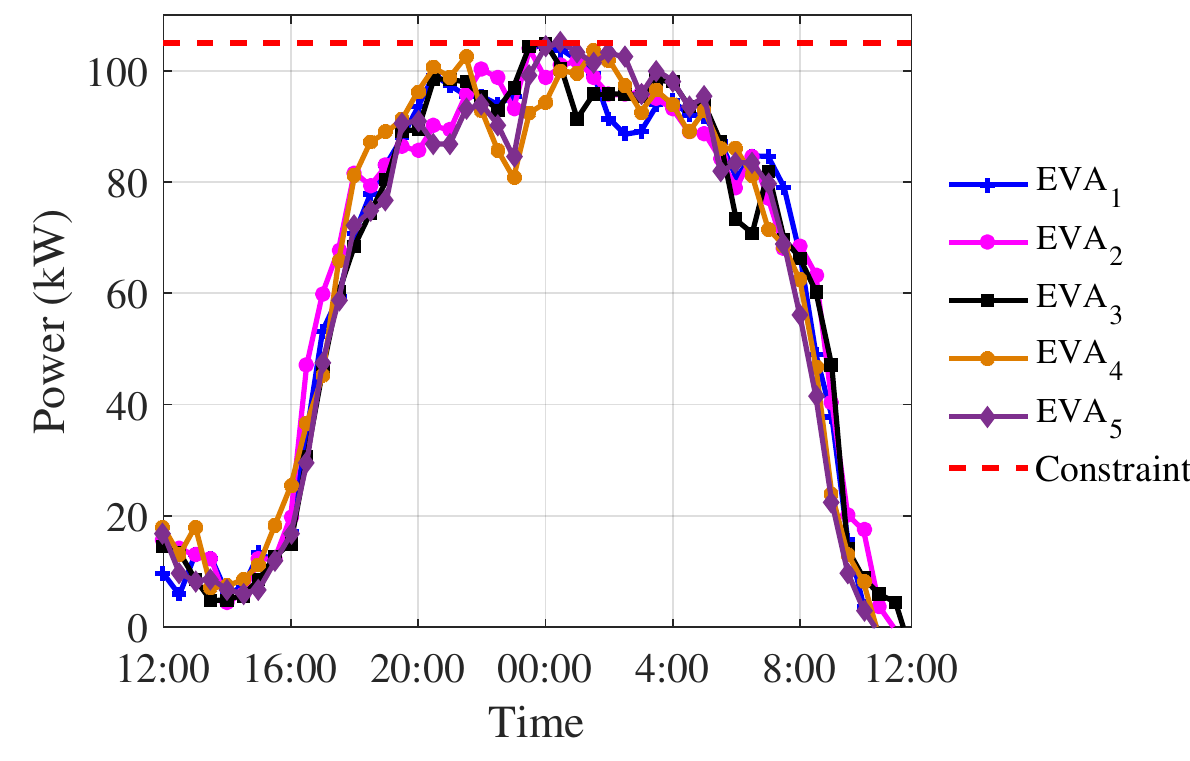}\label{CSat2}}
	\caption{EVAs' feeder constraint satisfaction in \textit{LVM+CR+BDR-1} (left) and \textit{LVM+CR+BDR-2} (right) modes by HDEVCS.}\label{CCompar2}
\end{figure}

As it is illustrated, the aggregated load profile in both modes is flattened by HDEVCS while the capacity constraints of EVAs' are not exceeded (\figref{CCompar2}). 

\begin{figure}[h!]
	\centering \includegraphics[clip,scale=0.45, trim=0cm 0cm 0cm 0cm] 
	{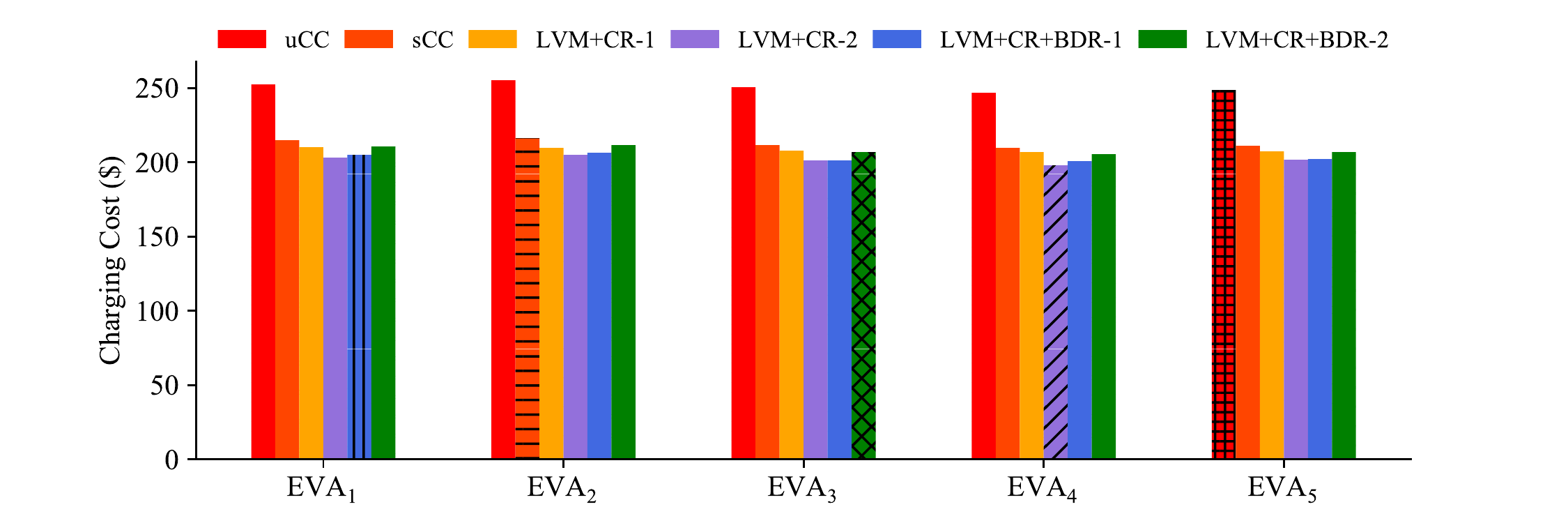}
	\caption{Aggregated EV charging cost of EVAs for different charging modes.}
	\label{CREVA}
\end{figure}

\begin{figure}[h!]
	\centering \includegraphics[clip,scale=0.45, trim=0cm 0cm 0cm 0cm] 
	{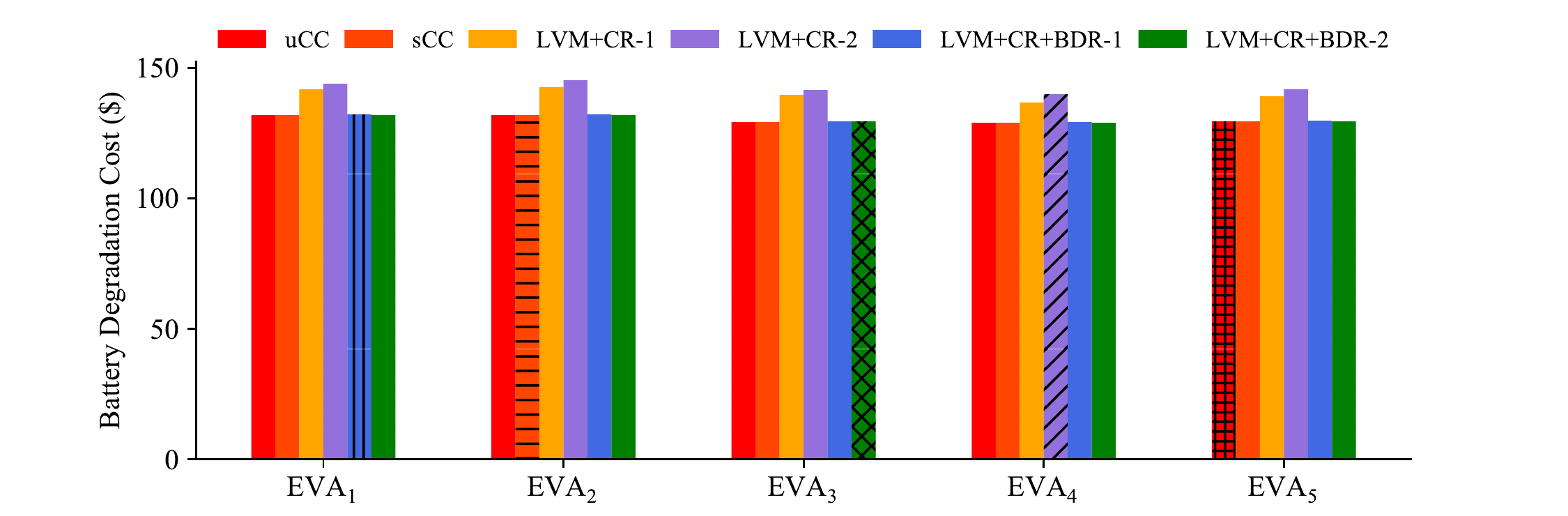}
	\caption{Aggregated EV battery degradation cost of EVAs for different charging modes.}
	\label{BDREVA}
\end{figure}

To compare the performance of HDEVCS in the first and second scenarios, the aggregated CR and BDR for EVAs is shown in \figref{CREVA} and \figref{BDREVA}, respectively. Since the weighting factor of CR objective function in \textit{LVM+CR-2} is larger than the weighting factor in \textit{LVM+CR-1}, the least aggregated CR for all EVAs is achieved by \textit{LVM+CR-2} (\figref{CREVA}).
As EVs start charging with the maximum power rating in \textit{uCC} mode, their charging cost is more than all other modes.
Considering \figref{BDREVA}, the least battery degradation cost is achieved by \textit{uCC}, \textit{sCC}, \textit{LVM+CR+BDR-1}, and \textit{LVM+CR+BDR-1}.
The reason is that the charging power profile in \textit{uCC} and \textit{sCC} is constant which according to \eqref{EVBObjS2} minimizes the battery degradation cost.
In \textit{LVM+CR+BDR-1} and \textit{LVM+CR+BDR-2} modes,  BDR cost is the objective function of EVs, therefore EVCS reduces BDR as well. However, battery degradation cost is larger in \textit{LVM+CR-1} and \textit{LVM+CR-2} as EVs try to reduce their charging cost only, and battery degradation reduction is not considered in the EVCS's optimization function. 

To further clarify the comparison of different simulated scenarios and the charging modes, the results are summarized in Table \ref{SimPerfMetrics} where the best obtained result for each metric is shown in bold.
As it is expected, \textit{LVM+CR-1} gives the lowest accumulated PTP, PTA, and RMS as the weighting factors of LVM and CR are equal, while the least aggregated charging cost is obtained by increasing the weight of CR in \textit{LVM+CR-2}.
As it was discussed, the lowest accumulated BDR is achieved by \textit{uCC}, \textit{sCC}, \textit{LVM+CR+BDR-1}, and \textit{LVM+CR+BDR-1}. The last column of the table shows the normalized accumulated objective values obtained by each charging mode. The weighting factors of the normalized values are equal to $1$. According to the calculated values, \textit{LVM-CR-1} charging mode results in the best accumulated performance, and \textit{uCC} charging mode leads to the worst performance.
  
\begin{figure}[h!]
	\centering \includegraphics[clip,scale=0.45, trim=1cm 9.2cm 0cm 10.6cm] 
	{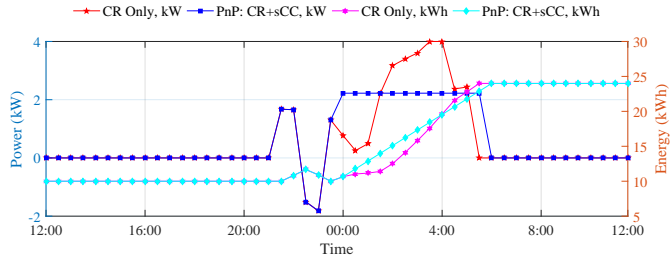}
	\caption{PnP using RH-HDEVCS: the EV agent switches from \textit{CR} mode to \textit{sCC} mode at $\texttt{23:30}$.}
	\label{fig:PnP}
\end{figure}

\begin{table}[h!]
	\caption{LVM Performance Metrics, CR, and BDR Improvement Using HDEVCS- \textit{System\textsubscript{1}}.}\label{SimPerfMetrics}
	\centering
	\footnotesize
	   \begin{threeparttable}[b]
	\setlength\tabcolsep{3pt}
	   \renewcommand{\arraystretch}{0.8}
	\begin{tabular}{@{} ccccccc @{}}
		\hlinewd{1pt}
		\addlinespace[0.05cm]
		EVCS Mode  & PTP(kW) & PTA & RMS(kW) & ACC\tnote{1} (\$) & BDC\tnote{2} (\$)& NAP\tnote{3}\\ \hlinewd{1pt}
		{uCC} & $1257$ & $0.21$ & $419.0$ & $1253$  & $\mathbf{651}$ & $4.91$  \\
		\addlinespace[0.05cm]
		{sCC} & $589$ & $0.10$ & $272.5$ & $1062$  &  $652$ & $3.34$ \\
		\addlinespace[0.05cm]
		{LVM+CR-1} & $\mathbf{494}$ & $\mathbf{0.08}$ & $\mathbf{258.2}$ & $1041$  & $699$ & $\mathbf{3.20}$   \\
		\addlinespace[0.05cm]
		{LVM+CR-2}& {$548$} & $0.09$ & $263.0$ &  $\mathbf{1008}$ & $712$ & $3.28$ \\
		\addlinespace[0.05cm]
		{LVM+CR+BDR-1}& {$549$} & $0.09$ & $263.6$ & $1015$  & $652$ & $3.21$ \\ 
		\addlinespace[0.05cm]
		{LVM+CR+BDR-2}& {$533$} & $0.09$ & $264.2$ & $1040$  & $\mathbf{651}$ & $3.22$   \\
		\addlinespace[0.05cm]
		\hlinewd{1pt}
	\end{tabular}
   \begin{tablenotes}
     \item[1] aggregated charging cost.
     \item[2] aggregated battery degradation cost.
     \item[3] normalized accumulated performance.
   \end{tablenotes}
  \end{threeparttable}
\end{table}
\paragraph{Plug-and-Play}
As mentioned before, the advantage of RH-HDEVCS is PnP in terms of the agents' objective function. That is, each agent may change its objective function in any RH iteration. This is illustrated in \figref{fig:PnP} where an EV is plugged in at $\texttt{21:00}$ when its battery energy is $10$ kWh, and it is unplugged at $\texttt{6:30}$ when it is fully charged. The EV is charged in CR mode until $\texttt{23:30}$ at which the EV owner switches the charging mode to \textit{sCC}, i.e. constant charging power. The desired energy stored in the battery at departure time is $24$ kWh. 

\subsection{\textit{System\textsubscript{2}}: Large-Scale Case Study}
\textit{System\textsubscript{2}} consists of $50$ EVAs each of which is supplying $180$ EVs.
The simulations are executed for two main scenarios with similar objective functions but different constraints.
The objective function in both scenarios includes LVM and CR. However, the feeder capacity constraints are considered in the first scenario (\textit{C/LVM+CR}), while there is no feeder constraint in the second scenario (\textit{UnC/LVM+CR}). In \textit{C/LVM+CR}, the maximum loading capacity is $180$ kW for the feeders supplying \textit{EVA}\textsubscript{$6$} and \textit{EVA}\textsubscript{$11$} and $175$ kW for the other feeders. The required charging energy, arrival time, and departure time of the EVs are generated as follows: the initial and designated EVs' battery energies are normally distributed over $[8,10]$ kWh and $[22,25]$ kWh, respectively; for $50\%$ of EVs, the arrival and departure times are normally distributed in $[\texttt{16:30},\texttt{20:30}]$ and $[\texttt{6:00},\texttt{9:30}]$, respectively; for the rest of EVs, the arrival and departure times are normally distributed in $[\texttt{6:00},\texttt{9:30}]$ and $[\texttt{16:30},\texttt{20:30}]$, respectively. 

\begin{figure}[h!]
	\centering \includegraphics[clip,scale=0.45, trim=1cm 9.2cm 0cm 10.0cm] 
	{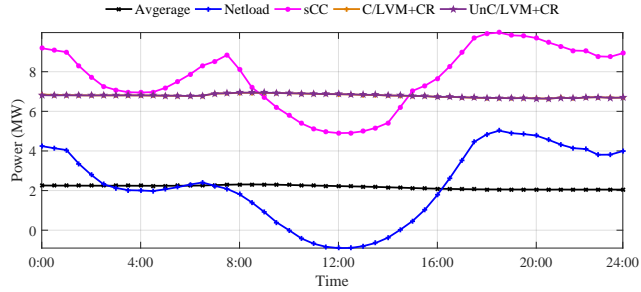}
	\caption{Aggregated load profile for \textit{sCC}, \textit{C/LVM+CR}, and \textit{UnC/LVM+CR}- \textit{System\textsubscript{2}}.}
	\label{fig:Scenario2LP}
\end{figure}

The aggregated load profiles obtained by \textit{sCC}, \textit{C/LVM+CR}, and \textit{UnC/LVM+CR} are shown in \figref{fig:Scenario2LP}. The load profiles of \textit{C/LVM+CR}, \textit{UnC/LVM+CR} coincides, and they perfectly fill in the valley (in $[\texttt{9:00},\texttt{15:00}]$) and shave the peak loads (at $\texttt{0:00}$, $\texttt{7:00}$ and $\texttt{18:30}$) which are seen in \textit{sCC} load profile.

\begin{figure}[h!]
	\centering \includegraphics[clip,scale=0.45, trim=0cm 9.2cm 1cm 10.4cm] 
	{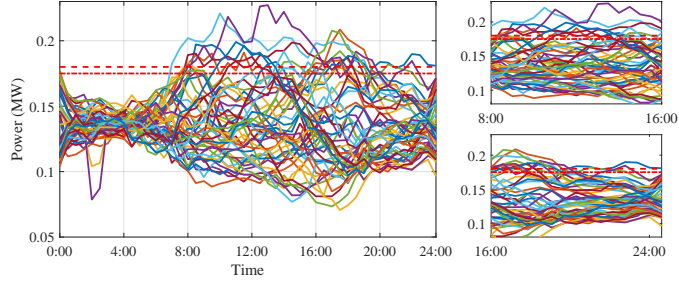}
	\caption{EVAs' feeder constraint violation in \textit{UnC/LVM+CR}.}
	\label{EVAP1}
\end{figure}

\begin{figure}[h!]
	\centering \includegraphics[clip,scale=0.45, trim=0cm 9.2cm 1cm 10.3cm] 
	{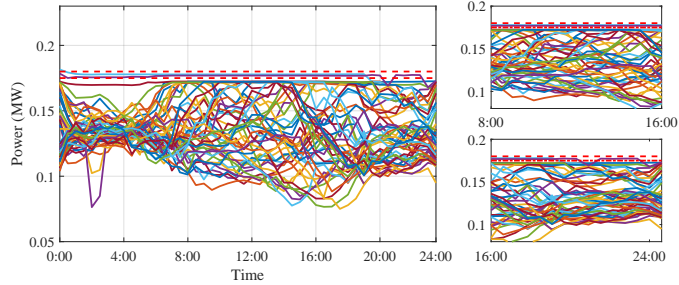}
	\caption{EVAs' feeder constraint satisfaction in \textit{C/LVM+CR}.}    
	\label{EVAP2}
\end{figure}

\begin{figure}[h!]
	\centering \includegraphics[clip,scale=0.45, trim=0cm 0cm 0cm 0cm] 
	{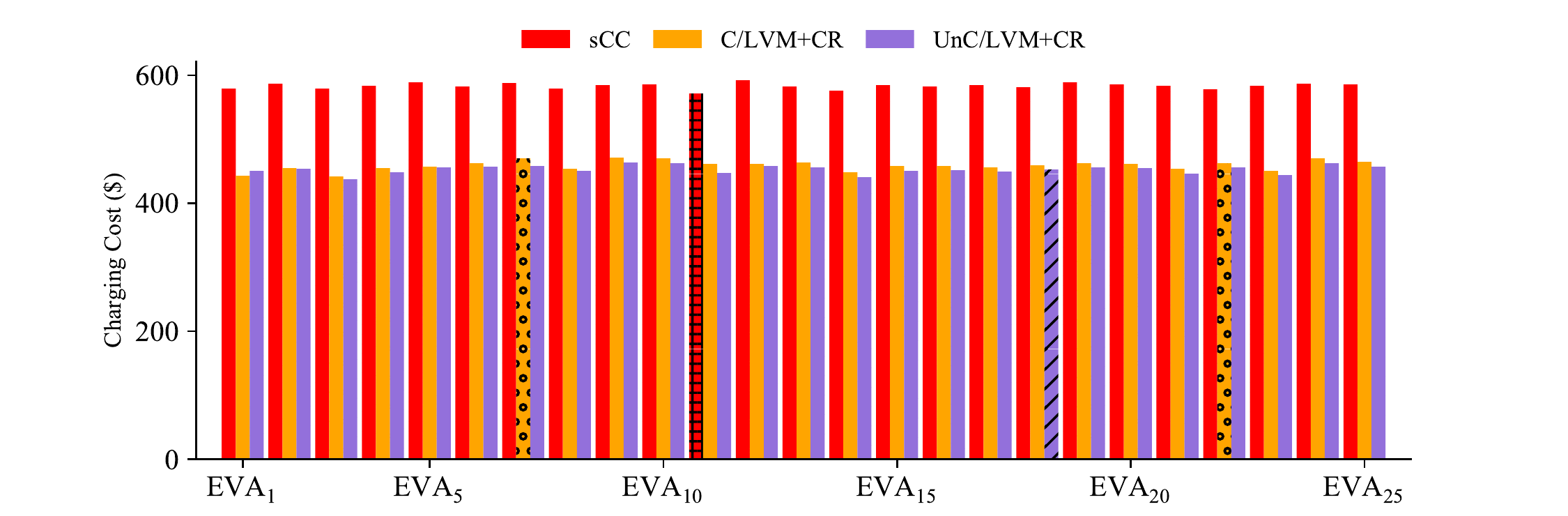}
	\caption{Aggregated charging cost of \textit{EVA}\textsubscript{$1$}-\textit{EVA}\textsubscript{$25$} by different charging modes for \textit{System\textsubscript{2}}.}
	\label{Price90001}
\end{figure}

\begin{figure}[h!]
	\centering \includegraphics[clip,scale=0.45, trim=0cm 0cm 0cm 0cm] 
	{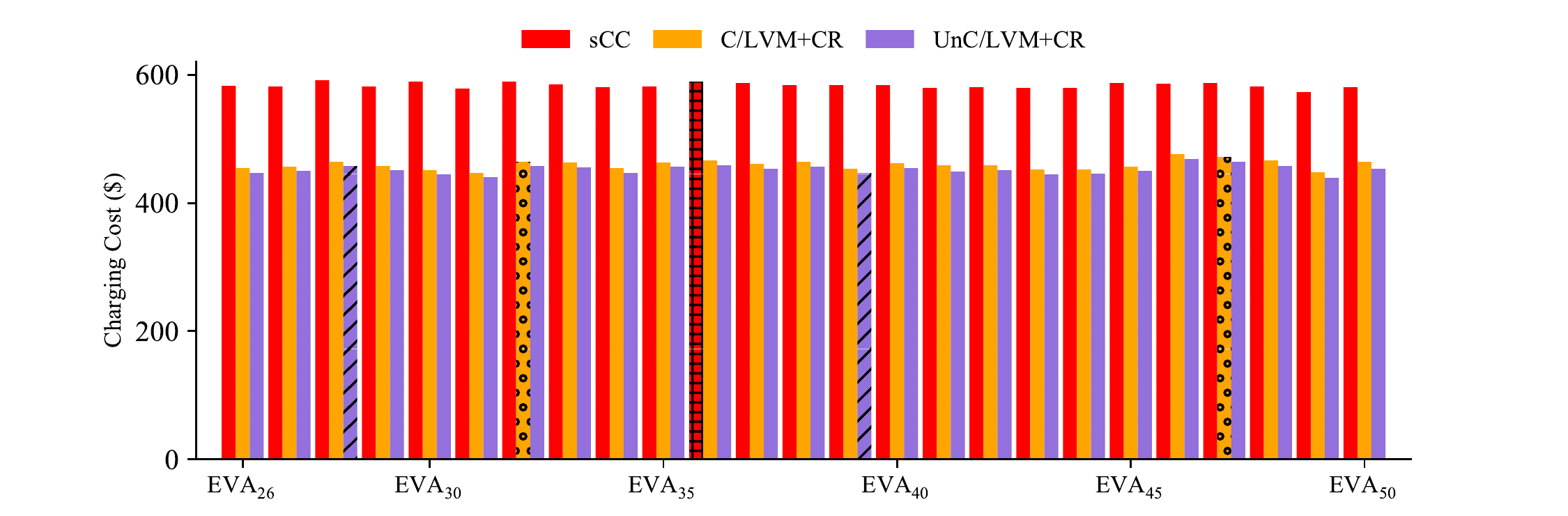}
	\caption{Aggregated charging cost of \textit{EVA}\textsubscript{$26$}-\textit{EVA}\textsubscript{$50$} by different charging modes for \textit{System\textsubscript{2}}.}
	\label{Price90002}
\end{figure}

Although the load profiles of \textit{C/LVM+CR} and \textit{UnC/LVM+CR} are similar, their EVAs' feeder load profiles are different. Comparing \figref{EVAP1} with \figref{EVAP2}, several EVA feeders are overloaded in \textit{UnC/LVM+CR} scenario while all the EVA feeder loads meet the constraints in \textit{C/LVM+CR}. This means that EVCS with feeder capacity constraints can accommodate a large population of EVs without any grid feeder expansion requirement. The other difference between \textit{C/LVM+CR} and \textit{UnC/LVM+CR} is seen by comparing the EVAs' aggregated charging costs in \figref{Price90001} and \figref{Price90002}. Although the difference is not considerable, \textit{UnC/LVM+CR} has less charging cost owing to the fact that EVs have more flexibility to shift their charging demand to the time periods with lower electricity price since there is no constraint on the EVAs' feeders.

The performance metrics calculated for the aggregated load profile in \textit{System\textsubscript{2}} are summarized in Table \ref{SimPerfMetrics2}, where ACC is calculated for the whole system. Comparing \textit{sCC} charging mode with \textit{C/LVM+CR}, improvement of the performance
metrics in this case is $94\%$, $30\%$ and $20\%$ for PTP, PTA, and RMS, respectively, and the cost reduction is $22\%$.
\begin{table}[h!]
	\caption{LVM Performance Metrics and CR Improvement Using HDEVCS- \textit{System\textsubscript{2}}.}\label{SimPerfMetrics2}
	\centering
	\footnotesize
	   \begin{threeparttable}[b]
	\setlength\tabcolsep{3pt}
	   \renewcommand{\arraystretch}{0.8}
	\begin{tabular}{@{} ccccc @{}}
		\hlinewd{1pt}
		\addlinespace[0.05cm]
		EVCS Mode  & PTP(kW) & PTA & RMS(kW) & ACC (\$)\\ \hlinewd{1pt}
		{sCC} & $5,073$ & $0.095$ & $5,747$ & $2,915$\\
		\addlinespace[0.05cm]
		{C/LVM+CR} & $292.4$ & $0.066$ & $4,608$ & $2,283$\\
		\addlinespace[0.05cm]
		{UnC/LVM+CR} & $301.3$ & $0.067$ & $4,609$ & $2,261$\\
		\addlinespace[0.05cm]
		\hlinewd{1pt}
	\end{tabular}
  \end{threeparttable}
\end{table}

\subsection{Comparison with Hierarchical ADMM \cite{VERSCHAE2016922}}
As it is already discussed, HDEVCS reduces the communication overhead and the convergence time compared to the methods in which the agents update their primal variable sequentially in two different steps. In this subsection, HDEVCS is compared with the hierarchical ADMM proposed by \citet{VERSCHAE2016922} which is designed based on the \textit{sharing problem} (\citet[Chapter~7.3]{Boyd1}). The hierarchical ADMM (\citet{VERSCHAE2016922}) consists of two layers of the sharing problem. The first layer is executed between DNO and EVAs, and the second layer between each EVA and its EVs. Therefore, the agents of ECI do not update their primal variable at the same time. To show the advantage of the proposed HDEVCS designed based on the \textit{exchange problem}, we run \textit{LVM+CR-1} charging mode for both \textit{System\textsubscript{1}} and \textit{System\textsubscript{2}} by the hierarchical ADMM. The convergence time and the number of iterations are shown in Table \ref{Comp12Layer}. For both methods, the penalty factor is similar ($\rho = 1$).
\begin{table}[t]
\caption{Comparison between HDEVCS and the hierarchical ADMM \cite{VERSCHAE2016922}.}\label{Comp12Layer}
\centering
\scalebox{0.80}{
\begin{tabular}{@{} cccc @{}}
\hlinewd{1pt}
 \textbf{Method} & \textbf{Case Study} &\textbf{Convergence Time (s)}  & \textbf{Iterations} \\ \hlinewd{1pt}
 \multirow{2}{*}{\textbf{HDEVCS}} & \textit{System\textsubscript{1}} &  $178.72$ & $64$ \\
& \textit{System\textsubscript{2}}  &  $3331.5$  &  $1193$ \\\hlinewd{0.5pt}
 \multirow{2}{*}{\textbf{Hierarchical ADMM \cite{VERSCHAE2016922}}} & \textit{System\textsubscript{1}}  & $446.8$ & $161$ \\
& \textit{System\textsubscript{2}}  &  $8204.1$  & $2938$ \\\hlinewd{1pt}
\end{tabular}}
\end{table}

As it is shown, HDEVCS improves convergence time by $60\%$ compared to the hierarchical ADMM. The convergence time, however, for \textit{System\textsubscript{2}} will be still considerable ($\approx 50$ min) if HDEVCS is embedded in RH as the netload dataset is collected every $30$ min, meaning that RH-HDEVCS should converge in less than $30$ min. However, it is worthy to mention that the computation times in Table \ref{Comp12Layer} are obtained by a CPU with $4$ cores, while RH-HDEVCS is proposed to be implemented in a multi-agent framework. Also, chances are the convergence time further improves using adaptive penalty term ($\rho$).

\section{Conclusion} \label{conclusion}
This paper has proposed a trilayer multi-agent framework for the optimal EV charging coordination to reduce the load variance and charging cost without violating the feeders' capacity constraints. By exploiting the configuration of the charging network and the mathematical properties of the EVCS problem, we have developed a novel hierarchical distributed method based on the \textit{exchange problem} for the optimal charging coordination problem that is solved by ADMM. Owing to the properties of the derived hierarchical \textit{exchange problem}, the second primal-update step of ADMM is eliminated, therefore all the agents update their primal variable in parallel, which results in the reduction of convergence time and iteration numbers. In addition, embedding the proposed method in the receding horizon feedback control gives flexibility to the agents to change their objective function in any receding horizon iteration, which is also called plug-and-play. To evaluate the performance of HDEVCS, it has been applied to two case studies, a small-scale and a large-scale system. The results have revealed that HDEVCS can reduce the peak load demand as well as EV charging and battery degradation costs significantly, while the grid feeders' capacity constraints are not violated. This means that the grid can accommodate a large population of EVs without investment in the grid capacity expansion. In comparison with the state-of-the-art trilayer charging framework, our proposed method reduces the convergence time and iteration numbers by $60\%$.  

The proposed HDEVCS, however, has some limitations which are considered in the authors' future works. The power flow model is not considered in the system model; although the load variance minimization improves voltage profile and reduces energy loss (\citet{Sortomme1}), the proposed method will need to be further modified if the minimum voltage is desired to be calculated precisely. Moreover, it is assumed that the communication latency is negligible in this work, which means the broadcast signals are received by the agents immediately. As in practice communication latency can be a bottleneck, specifically when the proposed algorithm is embedded in RH, the authors will further expand this work by modeling the communication delays to investigate if the algorithm can converge within a desired time. In that case, using ADMM with adaptive penalty factor can be also discussed. Another possible future direction of this research, recently investigated in the literature, is the battery swapping (\citet{Amiri1}) which is a solution for the battery charging time to facilitate EV penetration in the transportation system. Optimizing the performance of EV charging stations with this capability includes the quality of service to the customers and the charging cost. The authors will investigate the performance of HDEVCS by including the model of those charging stations with appropriate modification.    

\section{Appendix}\label{appendix}
The primal residual for each agent in $\text{\textbf{CL}}_{j}$ at iteration $k$ is obtained by
\begin{equation}  \label{primalres} 
\mathbf{r}_{j}^{k}=\overline{\mathbf{p}}_j^{k}.
\end{equation}
The dual residuals for $\textit{EV}_{v_i^j}$, $\textit{EVA}_{j}$ and DNO are calculated, respectively, by \eqref{dualres1}, \eqref{dualres2}-\eqref{dualres3} and \eqref{dualres4} as follows:
\begin{subequations} \label{dualres} 
\begin{align}
&\mathbf{s}_{v_i^j}^{k}=-\rho\mathcal{N}_{c_j}(\mathbf{p}_{v_i^j}^{k}-\mathbf{p}_{v_i^j}^{k-1}+\overline{\mathbf{p}}_j^{k-1}-\overline{\mathbf{p}}_j^{k}),~\forall i \in \mathbb{N}_v^j,~j \in \mathbb{N}_a \label{dualres1}\\
&\mathbf{s}_{au_j}^{k}=-\rho\mathcal{N}_{c_j}(-\mathbf{p}_{a_j}^{k}+\mathbf{p}_{a_j}^{k-1}+\overline{\mathbf{p}}_j^{k-1}-\overline{\mathbf{p}}_j^{k}),~\forall j \in \mathbb{N}_a\label{dualres2}\\
&\mathbf{s}_{a_j}^{k}=-\rho\mathcal{N}_{c}(\mathbf{p}_{a_j}^{k}-\mathbf{p}_{a_j}^{k-1}+\overline{\mathbf{p}}_{\mathcal{N}_c}^{k-1}-\overline{\mathbf{p}}_{\mathcal{N}_c}^{k}),~\forall j \in \mathbb{N}_a \label{dualres3}\\
&\mathbf{s}_{d}^{k}=-\rho\mathcal{N}_{c}(\mathbf{p}_{d}^{k}-\mathbf{p}_{d}^{k-1}+\overline{\mathbf{p}}_{\mathcal{N}_c}^{k-1}-\overline{\mathbf{p}}_{\mathcal{N}_c}^{k}). \label{dualres4}
\end{align}
\end{subequations}
\bibliography{mybibfile}

\end{document}